# A multirate Neumann-Neumann waveform relaxation method for heterogeneous coupled heat equations

Azahar Monge[*] and Philipp Birken[*]

[*]*Centre for Mathematical Sciences,
Lund University,
Box 118, 22100, Lund, Sweden*

## Abstract

An important challenge when coupling two different time dependent problems is to increase parallelization in time. We suggest a multirate Neumann-Neumann waveform relaxation algorithm to solve two heterogeneous coupled heat equations. In order to fix the mismatch produced by the multirate feature at the space-time interface a linear interpolation is constructed. The heat equations are discretized using a finite element method in space, whereas two alternative time integration methods are used: implicit Euler and SDIRK2. We perform a one-dimensional convergence analysis for the nonmultirate fully discretized heat equations using implicit Euler to find the optimal relaxation parameter in terms of the material coefficients, the stepsize and the mesh resolution. This gives a very efficient method which needs only two iterations. Numerical results confirm the analysis and show that the 1D nonmultirate optimal relaxation parameter is a very good estimator for the multirate 1D case and even for multirate and nonmultirate 2D examples using both implicit Euler and SDIRK2.

*Keywords:* Fluid-Structure Interaction, Coupled Problems, Transmission Problem, Domain Decomposition, Neumann-Neumann Method, Multirate, Thermal

# 1 Introduction

The main goal of this work is to describe a partitioned algorithm to solve two heterogeneous coupled heat equations allowing parallelization in time. In a partitioned approach different codes for the sub-problems are reused and the coupling is done by a master program which calls interface functions of the segregated codes [6, 7]. These algorithms are currently an active research topic driven by certain multiphysics applications where multiple physical models or multiple simultaneous physical phenomena involve solving coupled systems of partial differential equations (PDEs). An example of this is fluid structure interaction (FSI) [30, 4]. Moreover, we want that the time parallelization performed at the subsolvers works for different time grids. This is handled through multirate methods which are a classical field of research, see [10].

---

[*]e-mail: `azahar.monge@na.lu.se`; web page: http://www.maths.lu.se/staff/azahar-monge

Our prime motivation here is thermal interaction between fluids and structures, also called conjugate heat transfer. There are two domains with jumps in the material coefficients across the connecting interface. Conjugate heat transfer plays an important role in many applications and its simulation has proved essential [1]. Examples for thermal fluid structure interaction are cooling of gas-turbine blades, thermal anti-icing systems of airplanes [5], supersonic reentry of vehicles from space [24, 17], gas quenching, which is an industrial heat treatment of metal workpieces [16, 28] or the cooling of rocket nozzles [19, 20].

The classical way of parallelizing the numerical solution of PDEs is to use domain decomposition (DD) methods. These split the computational domain into subdomains and coordinate the coupling between the subdomains in an iterative manner. For an introduction to DD methods and their basic convergence results see [27, 29]. The Dirichlet-Neumann iteration is a standard DD method to find solutions of the coupled problem. The PDEs are solved sequentially using Dirichlet-, respectively Neumann boundary with data given from the solution of the other problem. Previous numerical experiments [2] showed that this iteration is fast for thermal FSI, and a convergence analysis of two heterogeneous linear heat equations showed that the fast behavior was a consequence of the strong jumps in the material coefficients [26].

In spite of the efficient behavior of the Dirichlet-Neumann iteration in the thermal FSI framework, it has two main disadvantages. Firstly, the subsolvers wait for each other, and therefore, they perform sequentially. Secondly, in the time dependent case the Dirichlet-Neumann iteration is used at each time step and consequently, both fields are solved with a common time resolution. Using instead a multirate scheme that allows for different time resolutions on each subdomain would be more efficient.

The aim of this work is to present a high order, parallel, multirate method for two heterogeneous coupled heat equations which could be applied to FSI problems. We use the Neumann-Neumann waveform relaxation (NNWR) method which is a waveform relaxation (WR) methods based on the classical Neumann-Neumann iteration [21, 13]. For time discretization we consider two alternatives, the implicit Euler method and a second order singly diagonally implicit Runge-Kutta method (SDIRK2). The WR methods were originally introduced by [22] for ordinary differential equation (ODE) systems, and they were used for the first time to solve time dependent PDEs in [14, 15]. They allow the use of different spatial and time discretizations for each subdomain which is specially useful in problems with strong jumps in the material coefficients [12] or the coupling of different models for the subdomains [11]. A time adaptive partitioned approach for thermal FSI was presented in [3]. In [23], two new iterative partitioned coupling methods that allow for the simultaneous execution of flow and structure solvers were introduced. However, parallelization in time for the coupling of heterogeneous materials was not considered.

Our algorithm has to take care of two aspects. On one hand, an interpolation procedure needs to be chosen to communicate data between the subdomains through the space-time interface in the multirate case. We want that the interpolation preserves a second order numerical solution of the coupled problem when using SDIRK2. On the other hand, the choice of the relaxation parameter for the NNWR method is crucial because when choosing the relaxation parameter right, two iterations are sufficient. In [21], a one-dimensional semidiscrete analysis shows that $\Theta = 1/4$ is the optimal relaxation parameter for two homogeneous coupled heat equations on two identical subdomains.

In this paper, we perform a fully discrete analysis of the NNWR algorithm for two heterogeneous coupled one-dimensional heat equations to find the optimal relaxation parameter in terms of the material coefficients. More specifically, we choose finite element methods (FEM) in space for both subdomains and implicit Euler method for the temporal discretization. Then, we derive the iteration matrix of the fully discrete NNWR algorithm with respect to the interface unknowns. In addition, we calculate the spectral radius of the iteration matrix through its eigendecomposition in order to estimate the optimal relaxation parameter $\Theta_{opt}$ which is dependent on the material coefficients, the time and space resolutions. In the case of homogeneous materials, $\Theta_{opt} = 1/4$ recovering the result in [21]. Furthermore, the asymptotic optimal relaxation parameters when approaching the continuous case in either time or space are also determined. In the spatial limit, the relaxation parameter turns out to be dependent on the heat conductivities, whereas in the temporal limit, we obtain dependency of the densities and the heat capacities.

In addition, we include numerical results where it is shown that the parallel, multirate method for two heterogeneous coupled heat equations introduced in this paper is extremely fast when choosing the right relaxation paramter. Moreover, we also show that the one-dimensional formula is a very good estimate for the multirate 1D case and even for multirate and nonmultirate 2D examples using both implicit Euler and SDIRK2. Finally, we also include a numerical comparison that shows that the NNWR method is a more efficient choice than the Dirichlet-Neumann waveform relaxation (DNWR) in the multirate case.

An outline of the paper now follows. In section 2, 3 and 4, we describe the model problem, the DNWR and the NNWR methods respectively. The FE space discretization is specified in section 5. In section 6, we describe the linear interpolation that needs to be performed at the space-time interface to get a multirate algorithm. Both time integration methods used in this paper are explained in section 7, these are implicit Euler and SDIRK2. In section 8, we present a derivation of the iteration matrix for a rather general discretization which is then applied to a specific 1D case in section 9. Numerical results are presented in section 10 and conclusions can be found in the last section.

## 2 Model problem

The unsteady transmission problem reads as follows, where we consider a domain $\Omega \subset \mathbb{R}^d$ which is cut into two subdomains $\Omega = \Omega_1 \cup \Omega_2$ with transmission conditions at the interface $\Gamma = \partial\Omega_1 \cap \partial\Omega_2$:

$$\begin{cases} \alpha_m \frac{\partial u_m(\mathbf{x},t)}{\partial t} - \nabla \cdot (\lambda_m \nabla u_m(\mathbf{x},t)) = 0, & \mathbf{x} \in \Omega_m \subset \mathbb{R}^d,\ m = 1,2, \\ u_m(\mathbf{x},t) = 0, & \mathbf{x} \in \partial\Omega_m \backslash \Gamma, \\ u_1(\mathbf{x},t) = u_2(\mathbf{x},t), & \mathbf{x} \in \Gamma, \\ \lambda_2 \frac{\partial u_2(\mathbf{x},t)}{\partial \mathbf{n}_2} = -\lambda_1 \frac{\partial u_1(\mathbf{x},t)}{\partial \mathbf{n}_1}, & \mathbf{x} \in \Gamma, \\ u_m(\mathbf{x},0) = u_m^0(\mathbf{x}), & \mathbf{x} \in \Omega_m, \end{cases} \quad (1)$$

where $t \in [T_0, T_f]$ and $\mathbf{n}_m$ is the outward normal to $\Omega_m$ for $m = 1, 2$.

The constants $\lambda_1$ and $\lambda_2$ describe the thermal conductivities of the materials on $\Omega_1$ and $\Omega_2$ respectively. $D_1$ and $D_2$ represent the thermal diffusivities of the materials and they are defined by

$$D_m = \frac{\lambda_m}{\alpha_m}, \quad \text{with} \quad \alpha_m = \rho_m c_{p_m} \quad (2)$$

where $\rho_m$ represents the density and $c_{p_m}$ the specific heat capacity of the material placed in $\Omega_m$, $m = 1, 2$.

## 3 The Dirichlet-Neumann Waveform Relaxation algorithm

The Dirichlet-Neumann waveform relaxation (DNWR) method is a basic iterative substructuring method in domain decomposition. The PDEs are solved sequentially using Dirichlet-, respectively Neumann boundary with data given from the solution of the other problem introduced in [13].

It starts with an initial guess $g^0(\mathbf{x}, t)$ on the interface $\Gamma \times (T_0, T_f]$, and then performs a three-step iteration. At each iteration $k$, imposing continuity of the solution across the interface, one first finds the local solution $u_1^{k+1}(\mathbf{x}, t)$ on $\Omega_1$ by solving the Dirichlet problem:

$$\begin{cases} \alpha_1 \frac{\partial u_1^{k+1}(\mathbf{x},t)}{\partial t} - \nabla \cdot (\lambda_1 \nabla u_1^{k+1}(\mathbf{x},t)) = 0, & \mathbf{x} \in \Omega_1, \\ u_1^{k+1}(\mathbf{x},t) = 0, & \mathbf{x} \in \partial\Omega_1 \backslash \Gamma, \\ u_1^{k+1}(\mathbf{x},t) = g^k(\mathbf{x},t), & \mathbf{x} \in \Gamma, \\ u_1^{k+1}(\mathbf{x},0) = u_1^0(\mathbf{x}), & \mathbf{x} \in \Omega_1. \end{cases} \quad (3)$$

Then, imposing continuity of the heat fluxes across the interface, one finds the local solution $u_2^{k+1}(\mathbf{x}, t)$ on $\Omega_2$ by solving the Neumann problem:

$$\begin{cases} \alpha_2 \frac{\partial u_2^{k+1}(\mathbf{x},t)}{\partial t} - \nabla \cdot (\lambda_2 \nabla u_2^{k+1}(\mathbf{x},t)) = 0, & \mathbf{x} \in \Omega_2, \\ u_2^{k+1}(\mathbf{x},t) = 0, & \mathbf{x} \in \partial\Omega_2 \backslash \Gamma, \\ \lambda_2 \frac{\partial u_2^{k+1}(\mathbf{x},t)}{\partial \mathbf{n_2}} = -\lambda_1 \frac{\partial u_1^{k+1}(\mathbf{x},t)}{\partial \mathbf{n_1}}, & \mathbf{x} \in \Gamma, \\ u_2^{k+1}(\mathbf{x},0) = u_2^0(\mathbf{x}), & \mathbf{x} \in \Omega_2. \end{cases} \quad (4)$$

Finally, the interface values are updated with

$$g^{k+1}(\mathbf{x},t) = \Theta u_2^{k+1}(\mathbf{x},t) + (1-\Theta)g^k(\mathbf{x},t), \quad \mathbf{x} \in \Gamma, \quad (5)$$

where $\Theta \in (0,1]$ is the relaxation parameter. The optimal relaxation parameter for the DNWR algorithm has been proved to be $\Theta = 1/2$ in [13] for the choice $\lambda_1 = \lambda_2 = \alpha_1 = \alpha_2 = 1$.

## 4 The Neumann-Neumann Waveform Relaxation algorithm

We now describe the Neumann-Neumann waveform relaxation (NNWR) algorithm [21]. The solution given by the NNWR method corresponds to the solution of the model problem (1) (proved in [18, chapt. 2]). The main advantage of the NNWR method is that it allows to find the solution on the subdomains in parallel.

The NNWR algorithm starts with an initial guess $g^0(\mathbf{x},t)$ on the space-time interface $\Gamma \times (T_0, T_f]$, and then performs a three-step iteration. At each iteration $k$, one first solves two Dirichlet problems on $\Omega_1$ and $\Omega_2$ simultaneously, then two Neumann problems are solved simultaneously again on $\Omega_1$ and $\Omega_2$, and finally, an update is performed to get a new guess $g^{k+1}(\mathbf{x},t)$ on the interface $\Gamma \times (T_0, T_f]$.

More specifically, imposing continuity of the solution across the interface (i.e, given a common initial guess $g^0(\mathbf{x},t)$ on $\Gamma \times (T_0, T_f)$), one can find the local solutions $u_m^{k+1}(\mathbf{x},t)$ on $\Omega_m$, $m = 1, 2$ through the following Dirichlet problems:

$$\begin{cases} \alpha_m \frac{\partial u_m^{k+1}(\mathbf{x},t)}{\partial t} - \nabla \cdot (\lambda_m \nabla u_m^{k+1}(\mathbf{x},t)) = 0, & \mathbf{x} \in \Omega_m, \\ u_m^{k+1}(\mathbf{x},t) = 0, & \mathbf{x} \in \partial\Omega_m \backslash \Gamma, \\ u_m^{k+1}(\mathbf{x},t) = g^k(\mathbf{x},t), & \mathbf{x} \in \Gamma, \\ u_m^{k+1}(\mathbf{x},0) = u_m^0(\mathbf{x}), & \mathbf{x} \in \Omega_m. \end{cases} \quad (6)$$

We now add into the framework the second coupling condition which is the continuity of the heat fluxes. To this end, one solves two simultaneous Neumann problems to get the correction functions $\psi_m^{k+1}(\mathbf{x},t)$ on $\Omega_m$, $m = 1, 2$ where the Neumann boundary condition at the interface $\Gamma \times (T_0, T_f)$ is prescribed by the

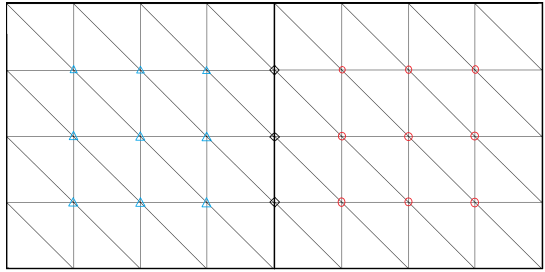

Figure 1: Splitting of $\Omega$ and finite element triangulation.

continuity of the heat fluxes of the solutions $u_m^{k+1}(\mathbf{x},t)$ given by the Dirichlet problems:

$$\begin{cases} \alpha_m \frac{\partial \psi_m^{k+1}(\mathbf{x},t)}{\partial t} - \nabla \cdot (\lambda_m \nabla \psi_m^{k+1}(\mathbf{x},t)) = 0, & \mathbf{x} \in \Omega_m, \\ \psi_m^{k+1}(\mathbf{x},t) = 0, & \mathbf{x} \in \partial\Omega_m \backslash \Gamma, \\ \lambda_m \frac{\partial \psi_m^{k+1}(\mathbf{x},t)}{\partial \mathbf{n}_m} = \lambda_1 \frac{\partial u_1^{k+1}(\mathbf{x},t)}{\partial \mathbf{n}_1} + \lambda_2 \frac{\partial u_2^{k+1}(\mathbf{x},t)}{\partial \mathbf{n}_2}, & \mathbf{x} \in \Gamma, \\ \psi_m^{k+1}(\mathbf{x},0) = 0, & \mathbf{x} \in \Omega_m. \end{cases} \quad (7)$$

Finally, the interface values are updated with

$$g^{k+1}(\mathbf{x},t) = g^k(\mathbf{x},t) - \Theta(\psi_1^{k+1}(\mathbf{x},t) + \psi_2^{k+1}(\mathbf{x},t)), \quad \mathbf{x} \in \Gamma, \quad (8)$$

where $\Theta \in (0,1]$ is the relaxation parameter. Note that choosing an appropriate relaxation parameter is crucial for the good performance of the NNWR algorithm [13]. If one uses the optimal relaxation parameter, two iterations are enough.

## 5 Semidiscrete method

We now describe a rather general space discretization of the problem (6)-(8). The core property we need is that the meshes of $\Omega_1$ and $\Omega_2$ share the same nodes on $\Gamma$ as shown in figure 1. Furthermore, we assume that there is a specific set of unknowns associated with the interface nodes. Otherwise, we allow at this point for arbitrary meshes on both sides.

Then, letting $\mathbf{u}_I^{(m)}, \psi_I^{(m)} : [T_0, T_f] \to \mathbb{R}^{S_m}$ where $S_m$ is the number of grid points on $\Omega_m$, $m=1,2$, and $\mathbf{u}_\Gamma, \psi_\Gamma^{(1)}, \psi_\Gamma^{(2)} : [T_0, T_f] \to \mathbb{R}^s$, where $s$ is the number of grid points at the interface $\Gamma$, we can write a general discretization of the first equation in (6) and (7), respectively, in a compact form as:

$$\mathbf{M}_{II}^{(m)} \dot{\mathbf{u}}_I^{(m),k+1}(t) + \mathbf{A}_{II}^{(m)} \mathbf{u}_I^{(m),k+1}(t) = -\mathbf{M}_{I\Gamma}^{(m)} \dot{\mathbf{u}}_\Gamma^k(t) - \mathbf{A}_{I\Gamma}^{(m)} \mathbf{u}_\Gamma^k(t), \quad (9)$$

$$\mathbf{M}_{II}^{(m)}\dot{\psi}_{I}^{(m),k+1}(t) + \mathbf{M}_{I\Gamma}^{(m)}\dot{\psi}_{\Gamma}^{(m),k+1}(t) + \mathbf{A}_{II}^{(m)}\psi_{I}^{(m),k+1}(t) + \mathbf{A}_{I\Gamma}^{(m)}\psi_{\Gamma}^{(m),k+1}(t) = \mathbf{0}, \tag{10}$$

where the initial conditions $\mathbf{u}_I^{(m)}(T_0), \psi_I^{(m)}(T_0) \in \mathbb{R}^{S_m}$ and $\mathbf{u}_\Gamma(T_0), \psi_\Gamma^{(m)}(T_0) \in \mathbb{R}^s$ for $m=1,2$ are known.

To close the system, we need an approximation of the normal derivatives on $\Gamma$. Letting $\phi_j$ be a nodal FE basis function on $\Omega_m$ for a node on $\Gamma$ we observe that the normal derivative of $u_m$ with respect to the interface can be written as a linear functional using Green's formula [29, pp. 3]. Thus, the approximation of the normal derivative is given by

$$\begin{aligned}\lambda_m \int_\Gamma \frac{\partial u_m}{\partial \mathbf{n}_m}\phi_j dS &= \lambda_m \int_{\Omega_m}(\Delta u_m \phi_j + \nabla u_m \nabla \phi_j)d\mathbf{x} \\ &= \alpha_m \int_{\Omega_m}\frac{d}{dt}u_m \phi_j + \lambda_m \int_{\Omega_m}\nabla u_m \nabla \phi_j d\mathbf{x}, \quad m=1,2.\end{aligned} \tag{11}$$

Consequently, the equation

$$\begin{aligned}&\mathbf{M}_{\Gamma\Gamma}^{(m)}\dot{\psi}_{\Gamma}^{(m),k+1}(t) + \mathbf{M}_{\Gamma I}^{(m)}\dot{\psi}_{I}^{(m),k+1}(t) + \mathbf{A}_{\Gamma\Gamma}^{(m)}\psi_{\Gamma}^{(m),k+1}(t) + \mathbf{A}_{\Gamma I}^{(m)}\psi_{I}^{(m),k+1}(t) \\ &= \sum_{i=1}^{2}\left(\mathbf{M}_{\Gamma\Gamma}^{(i)}\dot{\mathbf{u}}_\Gamma^k(t) + \mathbf{M}_{\Gamma I}^{(i)}\dot{\mathbf{u}}_I^{(i),k+1}(t) + \mathbf{A}_{\Gamma\Gamma}^{(i)}\mathbf{u}_\Gamma^k(t) + \mathbf{A}_{\Gamma I}^{(i)}\mathbf{u}_I^{(i),k+1}(t)\right), \quad m=1,2,\end{aligned} \tag{12}$$

is a discrete version of the third equation in (7) and completes the system (10).

We can now write a semidiscrete version of the NNWR algorithm using an ODE system. At each iteration $k$, one first solves the two Dirichlet problems in (9) obtaining $\mathbf{u}_I^{(m),k+1}(t)$ for $m=1,2$. Then, for the vector of unknowns $\psi_m^{k+1}(t) = \left(\psi_I^{(m),k+1}(t)^T \psi_\Gamma^{(m),k+1}(t)^T\right)^T$, one solves the following two Neumann problems in parallel that correspond to equations (10)-(12):

$$\mathbf{M}_m\dot{\psi}_m^{k+1}(t) + \mathbf{A}_m \psi_m^{k+1}(t) = \mathbf{b}^k, \quad m=1,2, \tag{13}$$

where

$$\mathbf{M}_m = \begin{pmatrix}\mathbf{M}_{II}^{(m)} & \mathbf{M}_{I\Gamma}^{(m)} \\ \mathbf{M}_{\Gamma I}^{(m)} & \mathbf{M}_{\Gamma\Gamma}^{(m)}\end{pmatrix}, \quad \mathbf{A}_m = \begin{pmatrix}\mathbf{A}_{II}^{(m)} & \mathbf{A}_{I\Gamma}^{(m)} \\ \mathbf{A}_{\Gamma I}^{(m)} & \mathbf{A}_{\Gamma\Gamma}^{(m)}\end{pmatrix}, \quad \mathbf{b}^k = \begin{pmatrix}\mathbf{0} \\ \mathbf{F}^k\end{pmatrix}, \tag{14}$$

with

$$\mathbf{F}^k = \sum_{i=1}^{2} \left( \mathbf{M}_{\Gamma\Gamma}^{(i)} \dot{\mathbf{u}}_\Gamma^k(t) + \mathbf{M}_{\Gamma I}^{(i)} \dot{\mathbf{u}}_I^{(i),k+1}(t) + \mathbf{A}_{\Gamma\Gamma}^{(i)} \mathbf{u}_\Gamma^k(t) + \mathbf{A}_{\Gamma I}^{(i)} \mathbf{u}_I^{(i),k+1}(t) \right). \quad (15)$$

Finally, the interfaces values are updated by

$$\mathbf{u}_\Gamma^{k+1}(t) = \mathbf{u}_\Gamma^k(t) - \Theta \left( \psi_\Gamma^{(1),k+1}(t) + \psi_\Gamma^{(2),k+1}(t) \right). \quad (16)$$

The iteration starts with some initial condition $\mathbf{u}_\Gamma^0(t)$ and a termination criterion must be chosen. One option would be $\|\mathbf{u}_\Gamma^{k+1}(t) - \mathbf{u}_\Gamma^k(t)\| \leq TOL$ where $TOL$ is a user defined tolerance. However, this option is memory consuming because it saves the solutions for all $t \in [T_0, T_f]$. Moreover, an extra interpolation step is needed in the multirate case, i.e, when having two nonconforming time grids. As we expect the error to be largest at the end point $T_f$ and because it simplifies the analysis to be presented for finding the optimal relaxation parameter, we propose the criterion $\|\mathbf{u}_\Gamma^{k+1}(T_f) - \mathbf{u}_\Gamma^k(T_f)\| \leq TOL$ where $T_f$ is the synchronization endpoint of the macrostep.

## 6 Space-time interface interpolation

The NNWR algorithm for parabolic problems was first introduced in [21, 13], but they do not consider the possibility of using two different step sizes on the two subdomains. In addition, their analysis does not include the coupling of two different materials. For those reasons, the goal of this paper is to introduce a parallel multirate method for the coupling of two heterogeneous heat equations and analyze its performance in the fully discrete case. This would be especially useful when coupling two different materials, where typically the field with higher heat conductivity needs a finer resolution than the other and therefore, efficiency will be gained by using a multirate method.

Note that both the Dirichlet and the Neumann problems in (9) and (13) allow the use of independent time discretization on each of the subdomains. Therefore, in the case of mismatched time grids, there exists the need to define an interface interpolation.

To this end, we consider a discrete problem in time with nonconforming time grids. Let $\tau_1 = \{t_1, t_2, .., t_{N_1}\}$ and $\tau_2 = \{t_1, t_2, .., t_{N_2}\}$ be two possibly different partitions of the time interval $[T_0, T_f]$ as shown in figure 2. We denote by $\Delta t_1 = (T_f - T_0)/N_1$ and $\Delta t_2 = (T_f - T_0)/N_2$ the two possibly different constant stepsizes corresponding to $\Omega_1$ and $\Omega_2$ respectively.

In order to exchange data at the space-time interface between the different time grids, we use a linear interpolation. For instance, if we want to interpolate the local discrete solution $G := (G_1, G_2, .., G_s) \in \mathbb{R}^{s \times N_2}$ from a given local discrete solution $F := (F_1, F_2, .., F_s) \in \mathbb{R}^{s \times N_1}$ at the space-time interface $\Gamma \times [T_0, T_f]$, with $s$ being the number of grid points at $\Gamma$, we use the following procedure: for each $k = 1, 2, .., s$ we consider the discrete vector $F_k := (F_{k,1}, F_{k,2}, .., F_{k,N_1}) \in \mathbb{R}^{N_1}$. Then, for each $t_i \in \tau_2$, $i = 1, .., N_2$, we find the subinterval $[t_j, t_{j+1}]$, $j = 0, .., N_1 - 1$ such that $t_i \in [t_j, t_{j+1}]$. We then define the linear interpolation polynomial $p_j(t)$ through the points $(t_j, F_{k,j})$ and $(t_{j+1}, F_{k,j+1})$, i.e, $p_j(t) = F_{k,j} + (t - t_j) \cdot (F_{k,j+1} - F_{k,j})/(t_{j+1} - t_j)$. Finally, we evaluate $p_j$ at $t_i$. Repeating the procedure for all the elements of $\tau_2$ we get the discrete vector $G_k$ whose components are given by,

$$G_{k,i} = \{p_j(t_i) \mid p_j : [t_j, t_{j+1}] \to \mathbb{R}, \forall t \in [t_j, t_{j+1}], j = 0, .., N_1 - 1\}, \quad (17)$$

for $i = 1, .., N_2$.

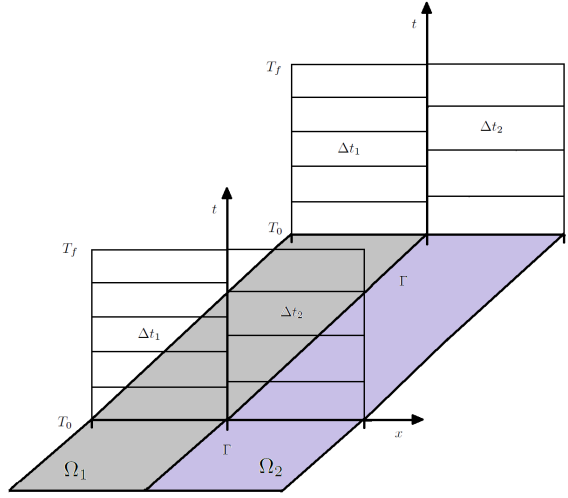

Figure 2: Nonconfoming time grids in the two-dimensional subdomains.

Thus, to interpolate $G$ from $F$ we use the interpolation function $G = I(\tau_2, \tau_1, F)$ explained above and summarized in algorithm 1. The same procedure can be applied to interpolate the discrete solution $F$ from a given discrete solution $G$.

**Algorithm 1** Interpolation to transfer data at the space-time interface.
1: **procedure** I($\tau_2, \tau_1, F$)
2:     **for** $k = 1, 2, .., s$ **do**
3:         **for** $t_i \in \tau_2$ **do**
4:             **for** $t_j \in \tau_1$ **do**
5:                 **if** $t_i \in [t_j, t_{j+1}]$ **then**
6:                       $G_{k,i} \leftarrow F_{k,j} + (t_i - t_j) \cdot (F_{k,j+1} - F_{k,j})/(t_{j+1} - t_j)$
        **return** $G$

# 7 Time integration

In this section we present a time discretized version of the NNWR method presented in equations (9), (13) and (16). In order to get a multirate algorithm we use a certain time integration method with time step $\Delta t_1 := (T_f - T_0)/N_1$ on $\Omega_1$ and with time step $\Delta t_2 := (T_f - T_0)/N_2$ on $\Omega_2$ and the interpolation presented in the previous section will be used to transfer data from one domain to the other. We let $n_m := 1, 2, .., N_m$ be the time integration indeces with respect to $\Omega_m$ and $t_{n_m}$ defines any time point of the grid for $m = 1, 2$. We have chosen two alternative time integration schemes as a basis to construct the multirate algorithm: the implicit Euler method and a second order singly diagonally implicit Runge-Kutta method (SDIRK2).

## 7.1 Implicit Euler

Applying the implicit Euler method with time step $\Delta t_1$ on $\Omega_1$ and with time step $\Delta t_2$ on $\Omega_2$ we can write the systems (9), (13) and (16) in a fully discrete form. At each fixed point iteration $k$, one first performs the time integration of the Dirichlet problems in parallel. Secondly, the interpolation explained in the previous section is used for the boundary conditions in order to solve the Neumann problems in parallel. Once the Neumann problems are solved, interpolation is again used to match the components of the update step. Finally, if the termination criterion is not fulfilled, one starts the process once more.

The local approximations and the solutions at the space-time interface are given by the vectors $\mathbf{u}_I^{(m),k,n_m} \approx \mathbf{u}_I^{(m),k}(t_{n_m}) \in \mathbb{R}^{S_m}$ and $\mathbf{u}_\Gamma^{k,n_m} \approx \mathbf{u}_\Gamma^k(t_{n_m}) \in \mathbb{R}^s$ respectively. Remember that $S_m$ is the number of spatial grid points on $\Omega_m$ and $s$ is the number of spatial grid points at the interface $\Gamma$. Similarly, the corrections both in the subdomains and at the interface are given by the vectors $\psi_I^{(m),k,n_m} \approx \psi_I^{(m),k}(t_{n_m}) \in \mathbb{R}^{S_m}$ and $\psi_\Gamma^{(m),k,n_m} \approx \psi_\Gamma^{(m),k}(t_{n_m}) \in \mathbb{R}^s$, respectively.

At each iteration $k$, one first solves the two Dirichlet problems from (9) for $n_m = 1, 2, .., N_m$, with $\mathbf{u}_I^{(m),k+1,0} \approx \mathbf{u}_I^{(m)}(T_0)$, $m = 1, 2$ and $\mathbf{u}_\Gamma^{k+1,0} \approx \mathbf{u}_\Gamma(T_0)$ simultaneously:

$$\left(\frac{\mathbf{M}_{II}^{(m)}}{\Delta t_m} + \mathbf{A}_{II}^{(m)}\right) \mathbf{u}_I^{(m),k+1,n_m+1} = -\left(\frac{\mathbf{M}_{I\Gamma}^{(m)}}{\Delta t_m} + \mathbf{A}_{I\Gamma}^{(m)}\right) \mathbf{u}_\Gamma^{k,n_m+1} \\ + \frac{\mathbf{M}_{II}^{(m)}}{\Delta t_m} \mathbf{u}_I^{(m),k+1,n_m} + \frac{\mathbf{M}_{I\Gamma}^{(m)}}{\Delta t_m} \mathbf{u}_\Gamma^{k,n_m}, \quad (18)$$

for $m = 1, 2$. Note that interpolation is not needed to solve the Dirichlet problems because $\mathbf{u}_I^{(1),k+1,n_1+1}$ in (18) is only dependent on terms related to $\Omega_1$. In the same way, $\mathbf{u}_I^{(2),k+1,n_2+1}$ in (18) only depends on $n_2$.

We compute now the fluxes $\tilde{\mathbf{F}}_1^{k,\tau_1} := \tilde{\mathbf{f}}_1^{k,\tau_1} + I(\tau_1, \tau_2, \tilde{\mathbf{f}}_2^{k,\tau_2})$ and $\tilde{\mathbf{F}}_2^{k,\tau_2} := \tilde{\mathbf{f}}_2^{k,\tau_2} + I(\tau_2, \tau_1, \tilde{\mathbf{f}}_1^{k,\tau_1})$ in (15) with

$$\tilde{\mathbf{f}}_m^{k,n_m} = \left(\frac{\mathbf{M}_{\Gamma\Gamma}^{(m)}}{\Delta t_m} + \mathbf{A}_{\Gamma\Gamma}^{(m)}\right)\mathbf{u}_\Gamma^{k,n_m+1} + \left(\frac{\mathbf{M}_{\Gamma I}^{(m)}}{\Delta t_m} + \mathbf{A}_{\Gamma I}^{(m)}\right)\mathbf{u}_I^{(m),k+1,n_m+1}$$
$$-\frac{\mathbf{M}_{\Gamma\Gamma}^{(m)}}{\Delta t_m}\mathbf{u}_\Gamma^{k,n_m} - \frac{\mathbf{M}_{\Gamma I}^{(m)}}{\Delta t_m}\mathbf{u}_I^{(m),k+1,n_m}, \quad (19)$$

where $n_m = 1,..,N_m$ and $\tau_m = \{t_1,t_2,..,t_{N_m}\}$ are the corresponding time grids on $\Omega_m$ for $m = 1,2$. Note that unlike in the Dirichlet problems, we need to use the interpolation described in the previous section. We use it to calculate $\tilde{\mathbf{F}}_1^{k,\tau_1}$ and $\tilde{\mathbf{F}}_2^{k,\tau_2}$ because their components run over different time integrations (one indicated by $n_1$ and the other by $n_2$).

One can now rewrite the Neumann problems in (13) in terms of the vector of unknowns $\psi_m^{k+1,n_m+1} := \left(\psi_I^{(m),k+1,n_m+1^T} \psi_\Gamma^{(m),k+1,n_m+1^T}\right)^T$. One then solves the two Neumann problems for $n_m = 1,2,..,N_m$, with $\psi_m^{k+1,0} \approx \psi_m(T_0)$, $m = 1,2$ in parallel:

$$\left(\frac{\mathbf{M}_m}{\Delta t_m} + \mathbf{A}_m\right)\psi_m^{k+1,n_m+1} = \frac{\mathbf{M}_m}{\Delta t_m}\psi_m^{k+1,n_m} + \tilde{\mathbf{b}}^{k,n_m}, \quad (20)$$

where $\tilde{\mathbf{b}}^{k,n_m} = \left(\mathbf{0}^T \tilde{\mathbf{F}}_m^{k,n_m T}\right)^T$.

Then, the interfaces values are updated respectively by

$$\mathbf{u}_\Gamma^{k+1,\tau_1} = \mathbf{u}_\Gamma^{k,\tau_1} - \Theta\left(\psi_\Gamma^{(1),k+1,\tau_1} + I\left(\tau_1,\tau_2,\psi_\Gamma^{(2),k+1,\tau_2}\right)\right), \quad (21)$$

$$\mathbf{u}_\Gamma^{k+1,\tau_2} = \mathbf{u}_\Gamma^{k,\tau_2} - \Theta\left(\psi_\Gamma^{(2),k+1,\tau_2} + I\left(\tau_2,\tau_1,\psi_\Gamma^{(1),k+1,\tau_1}\right)\right). \quad (22)$$

Here, interpolation is needed to perform the additions because $\psi_\Gamma^{(1),k+1,\tau_1}$ and $\psi_\Gamma^{(2),k+1,\tau_2}$ correspond to different time integrations.

Finally, if the termination criteria $\|\mathbf{u}_\Gamma^{k+1,N_m} - \mathbf{u}_\Gamma^{k,N_m}\| \approx \|\mathbf{u}_\Gamma^{k+1}(T_f) - \mathbf{u}_\Gamma^k(T_f)\|$ is not small enough, one starts the process from (18) once more.

Summarizing, figure 3 sketches the communication needed for the NNWR algorithm just explained. Algorithm 2 and 3 summarize the discrete Dirichlet solver in (18) and the discrete Neumann solver in (20) respectively. Furthermore, the complete NNWR algorithm is summarized in algorithm 4.

## 7.2 SDIRK2

As an alternative, we introduce here a higher order version of the multirate algorithm presented above. Specifically, we consider the second order singly

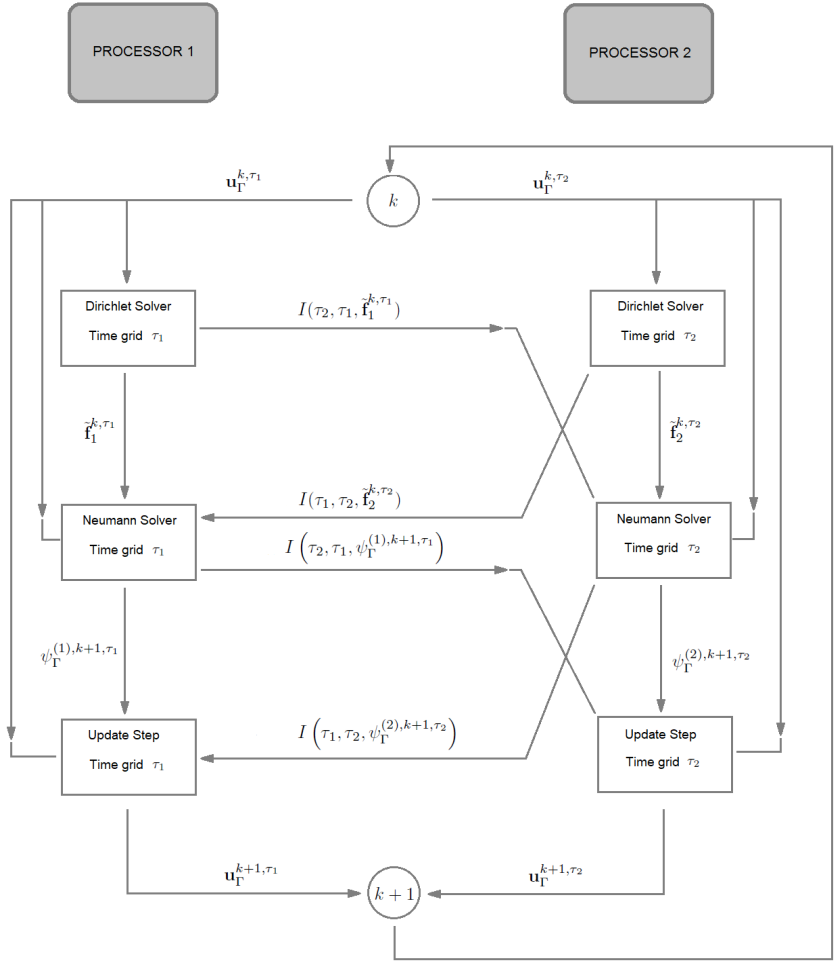

Figure 3: Illustration of the NNWR algorithm using implicit Euler. The process starts with the space-time interface functions $\mathbf{u}_\Gamma^{k,\tau_m}$, $\tau_m = \{t_1, t_2, .., t_{N_m}\}$ for $m = 1, 2$ corresponding to the two nonconforming time grids. Those are needed to run the Dirichlet solvers in parallel getting $\mathbf{u}_I^{(m),k+1,\tau_m}$, $m = 1, 2$. In order to run the Neumann solvers for the corrections of the solution, one needs to provide the fluxes $\tilde{\mathbf{f}}_1^{k,\tau_1}$, $\tilde{\mathbf{f}}_2^{k,\tau_2}$ and their corresponding interpolations $I(\tau_2, \tau_1, \tilde{\mathbf{f}}_1^{k,\tau_1})$, $I(\tau_1, \tau_2, \tilde{\mathbf{f}}_2^{k,\tau_2})$. One can then run the Neumann problems in parallel getting the corrections $\psi_\Gamma^{(1),k+1,\tau_1}$ and $\psi_\Gamma^{(2),k+1,\tau_2}$ at the space-time interface. Finally, those and their interpolations $I(\tau_1, \tau_2, \psi_\Gamma^{(2),k+1,\tau_2})$ and $I(\tau_2, \tau_1, \psi_\Gamma^{(1),k+1,\tau_1})$ are used to update the space-time interface values. If needed, the process is restarted.

**Algorithm 2** Solver for the Dirichlet problems in (18).

1: **procedure** SOLVEDIRICHLET($\mathbf{u}_I^{k+1,n}$, $\mathbf{u}_\Gamma^{k,n}$, $\mathbf{u}_\Gamma^{k,n+1}$)
2: $\quad \mathbf{u}_I^{k+1,n+1} \leftarrow$ SOLVELINEARSYSTEM($\mathbf{u}_I^{k+1,n}$, $\mathbf{u}_\Gamma^{k,n}$, $\mathbf{u}_\Gamma^{k,n+1}$), # *solve* (18)
3: $\quad \tilde{\mathbf{f}}^{k,n} \leftarrow \left(\frac{\mathbf{M}_{\Gamma\Gamma}}{\Delta t} + \mathbf{A}_{\Gamma\Gamma}\right) \mathbf{u}_\Gamma^{k,n+1} + \left(\frac{\mathbf{M}_{\Gamma I}}{\Delta t} + \mathbf{A}_{\Gamma I}\right) \mathbf{u}_I^{k+1,n+1} - \frac{\mathbf{M}_{\Gamma\Gamma}}{\Delta t} \mathbf{u}_\Gamma^{k,n} - \frac{\mathbf{M}_{\Gamma I}}{\Delta t} \mathbf{u}_I^{k+1,n}$, # *compute* (19)
$\quad\quad$ **return** $\mathbf{u}_I^{k+1,n+1}, \tilde{\mathbf{f}}^{k,n}$

---

**Algorithm 3** Solver for the Neumann problems in (20).

1: **procedure** SOLVENEUMANN($\psi_I^{k+1,n}$, $\psi_\Gamma^{k+1,n}$, $\tilde{\mathbf{F}}^{k,n}$)
2: $\quad \psi_I^{k+1,n+1}, \psi_\Gamma^{k+1,n+1} \leftarrow$ SOLVELINEARSYSTEM($\psi_I^{k+1,n}, \psi_\Gamma^{k+1,n}, \tilde{\mathbf{F}}^{k,n}$), # *solve* (20)
$\quad\quad$ **return** $\psi_I^{k+1,n+1}, \psi_\Gamma^{k+1,n+1}$

---

**Algorithm 4** NNWR algorithm using implicit Euler.

1: **procedure** NNWR($\tau_1$, $\tau_2$, $\alpha_1$, $\alpha_2$, $\lambda_1$, $\lambda_2$, $\Theta$, $TOL$)
2: $\quad \mathbf{u}_I^{(m),k+1,0}, \mathbf{u}_\Gamma^0(\tau_m), \psi_I^{(m),k+1,n_m+1}, \psi_\Gamma^{(m),k+1,n_m+1} \leftarrow$ INITIALIZATION
3: $\quad$ **while** $\|\mathbf{u}_\Gamma^{k+1,N_m} - \mathbf{u}_\Gamma^{k,N_m}\| \leq TOL$ **do**
4: $\quad\quad$ **for** $t_{n_1} \in \tau_1$ **do**
5: $\quad\quad\quad \mathbf{u}_I^{(1),k+1,n_1+1}, \tilde{\mathbf{f}}_1^{k,n_1} \leftarrow$ SOLVEDIRICHLET($\mathbf{u}_I^{(1),k+1,n_1}, \mathbf{u}_\Gamma^{k,n_1}, \mathbf{u}_\Gamma^{k,n_1+1}$)
6: $\quad\quad$ **for** $t_{n_2} \in \tau_2$ **do** (in parallel to 4)
7: $\quad\quad\quad \mathbf{u}_I^{(2),k+1,n_2+1}, \tilde{\mathbf{f}}_2^{k,n_2} \leftarrow$ SOLVEDIRICHLET($\mathbf{u}_I^{(2),k+1,n_2}, \mathbf{u}_\Gamma^{k,n_2}, \mathbf{u}_\Gamma^{k,n_2+1}$)
8: $\quad\quad \tilde{\mathbf{F}}_1^{k,\tau_1} \leftarrow \tilde{\mathbf{f}}_1^{k,\tau_1} + \mathrm{I}(\tau_1, \tau_2, \tilde{\mathbf{f}}_2^{k,\tau_2})$
9: $\quad\quad \tilde{\mathbf{F}}_2^{k,\tau_2} \leftarrow \tilde{\mathbf{f}}_2^{k,\tau_2} + \mathrm{I}(\tau_2, \tau_1, \tilde{\mathbf{f}}_1^{k,\tau_1})$ (in parallel to 8)
10: $\quad\quad$ **for** $t_{n_1} \in \tau_1$ **do**
11: $\quad\quad\quad \psi_I^{(1),k+1,n_1+1}, \psi_\Gamma^{(1),k+1,n_1+1} \leftarrow$ SOLVENEUMANN($\psi_I^{(1),k+1,n_1}, \psi_\Gamma^{(1),k+1,n_1}, \tilde{\mathbf{F}}_1^{k,n_1}$)
12: $\quad\quad$ **for** $t_{n_2} \in \tau_2$ **do** (in parallel to 10)
13: $\quad\quad\quad \psi_I^{(2),k+1,n_2+1}, \psi_\Gamma^{(2),k+1,n_2+1} \leftarrow$ SOLVENEUMANN($\psi_I^{(2),k+1,n_2}, \psi_\Gamma^{(2),k+1,n_2}, \tilde{\mathbf{F}}_2^{k,n_2}$)
14: $\quad\quad \mathbf{u}_\Gamma^{k+1,\tau_1} \leftarrow \mathbf{u}_\Gamma^{k,\tau_1} - \Theta\left(\psi_\Gamma^{(1),k+1,\tau_1} + \mathrm{I}(\tau_1, \tau_2, \psi_\Gamma^{(2),k+1,\tau_2})\right)$
15: $\quad\quad \mathbf{u}_\Gamma^{k+1,\tau_2} \leftarrow \mathbf{u}_\Gamma^{k,\tau_2} - \Theta\left(\psi_\Gamma^{(2),k+1,\tau_2} + \mathrm{I}(\tau_2, \tau_1, \psi_\Gamma^{(1),k+1,\tau_1})\right)$ (in p. to 14)

diagonally implicit Runge-Kutta (SDIRK2) as a basis to discretize the systems (9), (13) and (16) in time. Consider an autonomous initial value problem

$$\dot{\mathbf{u}}(t) = \mathbf{f}(\mathbf{u}(t)), \quad \mathbf{u}(0) = \mathbf{u}_0. \tag{23}$$

An SDIRK method is then defined as

$$\mathbf{U}^i = \mathbf{u}^n + \Delta t_n \sum_{k=1}^{i} a_{ik}\mathbf{f}(\mathbf{U}^k), \quad i = 1, .., j$$

$$\mathbf{u}^{n+1} = \mathbf{u}^n + \Delta t_n \sum_{i=1}^{j} b_i \mathbf{f}(\mathbf{U}^i) \tag{24}$$

with given coefficients $a_{ik}$ and $b_i$. The two-stage method SDIRK2 is defined by the coefficients in the following Butcher array:

$$\begin{array}{c|cc} a & a & 0 \\ 1 & 1-a & a \\ \hline & 1-a & a \end{array}$$

with $a = 1 - \frac{1}{2}\sqrt{2}$. As the coefficients $a_{2i}$ and $b_i$ for $i = 1, 2$ are identical, the second equation in (24) is superfluous because $\mathbf{u}^{n+1} = \mathbf{U}^2$.

The vectors $\mathbf{k}_i = \mathbf{f}(\mathbf{U}^i)$ are called stage derivatives and $j$ is the number of stages. Since the starting vector

$$\mathbf{s}_i = \mathbf{u}^n + \Delta t_n \sum_{k=1}^{i-1} a_{ik}\mathbf{k}_k, \quad i = 1, .., j-1, \tag{25}$$

is known, (24) is just a sequence of implicit Euler steps.

Applying SDIRK2 with time step $\Delta t_1$ on $\Omega_1$ and with time step $\Delta t_2$ on $\Omega_2$ we can write the systems (9), (13) and (16) in a fully discrete form. This algorithm preserves more or less the same structure as the one presented above for implicit Euler. The main difference lies in the fact that now both the Dirichlet and the Neumann solvers have to take into account the two stages of SDIRK2 as well as the interpolation has to be applied for each stage.

Therefore, at each fixed point iteration $k$, let $\mathbf{s}_1^{(m)} = \mathbf{u}_I^{(m),k+1,n_m}$ and $\mathbf{s}_2^{(m)} = \mathbf{u}_I^{(m),k+1,n_m} + \Delta t_m(1-a)\mathbf{k}_1^{(m)}$ be the starting vectors. Then, one first solves the two Dirichlet problems for $n_m = 1, 2, .., N_m$, with $\mathbf{u}_I^{(m),k+1,0}$, $m = 1, 2$, $\mathbf{u}_\Gamma^{k+1,0}$ simultaneously:

$$\left(\frac{\mathbf{M}_{II}^{(m)}}{a\Delta t_m} + \mathbf{A}_{II}^{(m)}\right)\mathbf{U}_j^{(m)} = \frac{\mathbf{M}_{II}^{(m)}}{a\Delta t_m}\mathbf{s}_j^{(m)} - \mathbf{M}_{I\Gamma}^{(m)}\dot{\mathbf{u}}_\Gamma^{k,n_m+j-1+(2-j)a}$$
$$- \mathbf{A}_{I\Gamma}^{(m)}\mathbf{u}_\Gamma^{k,n_m+j-1+(2-j)a}, \quad j = 1, 2. \tag{26}$$
$$\mathbf{u}_I^{(m),k+1,n_m+1} = \mathbf{U}_2^{(m)},$$

where $\mathbf{u}_I^{(m),k,n_m}$, $\mathbf{U}_j^{(m)}$, $\mathbf{s}_j^{(m)}$, $\mathbf{k}_j^{(m)} \in \mathbb{R}^{S_m}$ and $\mathbf{u}_\Gamma^{k,n_m} \in \mathbb{R}^s$. The stage derivatives are given by $\mathbf{k}_j^{(m)} = \frac{1}{a\Delta t_m}(\mathbf{U}_j^{(m)} - \mathbf{s}_j^{(m)})$. Note that the index $m = 1, 2$, denotes the subdomain and the index $j = 1, 2$, denotes the stage.

We compute now the fluxes $\mathbf{F}_j^{(1),k,\tau_1} := \mathbf{f}_j^{(1),k,\tau_1} + I(\tau_1, \tau_2, \mathbf{f}_j^{(2),k,\tau_2})$, $\mathbf{F}_j^{(2),k,\tau_2} := \mathbf{f}_j^{(2),k,\tau_2} + I(\tau_2, \tau_1, \mathbf{f}_j^{(1),k,\tau_1})$ in (15) with

$$\begin{aligned} \mathbf{f}_j^{(m),k,n_m} = &\mathbf{M}_{\Gamma\Gamma}^{(m)} \dot{\mathbf{u}}_\Gamma^{k,n_m+j-1+(2-j)a} + \mathbf{M}_{\Gamma I}^{(m)} \mathbf{k}_j^{(m)} \\ &+ \mathbf{A}_{\Gamma\Gamma}^{(m)} \mathbf{u}_\Gamma^{k,n_m+j-1+(2-j)a} + \mathbf{A}_{\Gamma I}^{(m)} \mathbf{U}_j^{(m)}, \end{aligned} \quad (27)$$

for $m = 1, 2$. Note that interpolation here is needed because the components of $\mathbf{F}_j^{(1),k,\tau_1}$ and $\mathbf{F}_j^{(2),k,\tau_2}$ for the two stages $j = 1, 2$ correspond to different time integrations.

One can now rewrite the Neumann problems in (13) in terms of the vector of unknowns $\psi_m^{k+1,n_m+1} := \left(\psi_I^{(m),k+1,n_m+1\,T} \psi_\Gamma^{(m),k+1,n_m+1\,T}\right)^T$ where $\psi_I^{(m),k+1,n_m+1} \in \mathbb{R}^{S_m}$ and $\psi_\Gamma^{(m),k+1,n_m+1} \in \mathbb{R}^s$. Let $\mathbf{s}_1^{(m)} = \psi_m^{k+1,n_m}$ and $\mathbf{s}_2^{(m)} = \psi_m^{k+1,n_m} + \Delta t_m (1-a) \mathbf{k}_1^{(m)}$ be the starting vectors. One then solves the two Neumann problems for $n_m = 1, 2, .., N_m$, with $\psi_m^{k+1,0} = \psi_m^{k+1}(T_0)$, $m = 1, 2$ in parallel:

$$\begin{aligned} \left(\frac{\mathbf{M}_m}{a\Delta t_m} + \mathbf{A}_m\right) \mathbf{Y}_j^{(m)} &= \frac{\mathbf{M}_m}{a\Delta t_m} \mathbf{s}_j^{(m)} + \mathbf{b}_j^{(m),k,n_m}, \quad j = 1, 2, \\ \psi_m^{k+1,n_m+1} &= \mathbf{Y}_2^{(m)}, \end{aligned} \quad (28)$$

where $\mathbf{Y}_j^{(m)}$, $\mathbf{s}_j^{(m)}$, $\mathbf{b}_j^{(m),k,n_m}$, $\mathbf{k}_j^{(m)} \in \mathbb{R}^{S_m+s}$, $\mathbf{k}_j^{(m)} = \frac{1}{a\Delta t_m}(\mathbf{Y}_j^{(m)} - \mathbf{s}_j^{(m)})$ and $\mathbf{b}_j^{(m),k,n_m} = \left(\mathbf{0}^T \ \mathbf{F}_j^{(m),k,n_m\,T}\right)^T$.

Then, the interfaces values are updated respectively by

$$\mathbf{u}_\Gamma^{k+1,\tau_1} = \mathbf{u}_\Gamma^{k,\tau_1} - \Theta\left(\psi_\Gamma^{(1),k+1,\tau_1} + I\left(\tau_1, \tau_2, \psi_\Gamma^{(2),k+1,\tau_2}\right)\right), \quad (29)$$

$$\mathbf{u}_\Gamma^{k+1,\tau_2} = \mathbf{u}_\Gamma^{k,\tau_2} - \Theta\left(\psi_\Gamma^{(2),k+1,\tau_2} + I\left(\tau_2, \tau_1, \psi_\Gamma^{(1),k+1,\tau_1}\right)\right). \quad (30)$$

Here, interpolation is needed because $\psi_\Gamma^{(1),k+1,\tau_1}$ and $\psi_\Gamma^{(2),k+1,\tau_2}$ are nonconforming.

Finally, if the termination criteria $\|\mathbf{u}_\Gamma^{k+1,N_m} - \mathbf{u}_\Gamma^{k,N_m}\|$ is not small enough, one starts the process from (26) once more.

We use a linear interpolation through the points $(t_{n_m}, \mathbf{u}_\Gamma^{k,n_m})$ and $(t_{n_m} + \Delta t_m, \mathbf{u}_\Gamma^{k,n_m+1})$ in order to approximate $\mathbf{u}_\Gamma^{k,n_m+a}$ in the first equation of (26) and in the fluxes (27), i.e:

$$\mathbf{u}_\Gamma^{k,n_m+a} \approx \mathbf{u}_\Gamma^{k,n_m} + a\left(\mathbf{u}_\Gamma^{k,n_m+1} - \mathbf{u}_\Gamma^{k,n_m}\right). \tag{31}$$

Furthermore, there are first order time derivatives in the first equation of (26) and in (27). We use forward differences to approximate all the remaining first order derivatives:

$$\dot{\mathbf{u}}_\Gamma^{k,n_m+j-1+(2-j)a} \approx \frac{\mathbf{u}_\Gamma^{k,n_m+1} - \mathbf{u}_\Gamma^{k,n_m}}{\Delta t_m}, \tag{32}$$

for $j = 1, 2$ and $m = 1, 2$.

Summarizing, the SDIRK2-NNWR algorithm just presented has the same structure as the implicit Euler NNWR algorithm described previously and sketched in figure 3. The main difference is that the whole procedure is repeated twice, once for each stage.

Algorithm 5 and 6 summarize the discrete Dirichlet solver in (26) and the discrete Neumann solver in (28) respectively. Furthermore, the complete SDIRK2-NNWR algorithm is summarized in algorithm 7.

---

**Algorithm 5** Solver for the Dirichlet problems in (26).

1: **procedure** SDIRK2DIRICHLET($\mathbf{u}_I^{k+1,n}$, $\mathbf{u}_\Gamma^{k,n}$, $\mathbf{u}_\Gamma^{k,n+1}$)
2:      $\mathbf{u}_\Gamma^{k,n+a} \leftarrow \mathbf{u}_\Gamma^{k,n} + a\left(\mathbf{u}_\Gamma^{k,n+1} - \mathbf{u}_\Gamma^{k,n}\right)$
3:      $\dot{\mathbf{u}}_\Gamma^{k,n+a}, \dot{\mathbf{u}}_\Gamma^{k,n+1} \leftarrow \left(\mathbf{u}_\Gamma^{k,n+1} - \mathbf{u}_\Gamma^{k,n}\right)/\Delta t$
4:      **for** $j = 1, 2$ **do**   # loop over stages
5:          $\mathbf{s}_j \leftarrow \mathbf{u}_I^{k+1,n} + \Delta t \sum_{l=1}^{j-1}(1-a)\mathbf{k}_l$
6:          $\mathbf{U}_j \leftarrow$ SOLVELINEARSYSTEM($\mathbf{s}_j, \dot{\mathbf{u}}_\Gamma^{k,n+j-1+(2-j)a}, \mathbf{u}_\Gamma^{k,n+j-1+(2-j)a}$),   # solve 1st eq in (26)
7:          $\mathbf{k}_j \leftarrow \frac{1}{a\Delta t}(\mathbf{U}_j - \mathbf{s}_j)$
8:      $\mathbf{u}_I^{k+1,n+1} \leftarrow \mathbf{U}_2$
9:      **for** $j = 1, 2$ **do**   # compute fluxes in (27)
10:         $\mathbf{f}_j^{k,n} \leftarrow \mathbf{M}_{\Gamma\Gamma}\dot{\mathbf{u}}_\Gamma^{k,n+j-1+(2-j)a} + \mathbf{M}_{\Gamma I}\mathbf{k}_j + \mathbf{A}_{\Gamma\Gamma}\mathbf{u}_\Gamma^{k,n+j-1+(2-j)a} + \mathbf{A}_{\Gamma I}\mathbf{U}_j$
     **return** $\mathbf{u}_I^{k+1,n+1}$, $\mathbf{f}_j^{k,n}$

---

## 8 Derivation of the iteration matrix

We are interested in the performance of the NNWR algorithm. As the rate of convergence of a linear iteration is given by the spectral radius of its iteration matrix, we derive in this section the iteration matrix with respect to the set of unknowns at the space-time interface for implicit Euler. A similar analysis to find the convergence rates of the Dirichlet-Neumann iteration for the unsteady transmission problem can be found in [26]. We intentionally avoid a derivation for SDIRK2 and we will show in the numerical results section that NNWR-SDIRK2 behaves as predicted by the analysis of implicit Euler. From now on we assume that we have conforming time grids, i.e, $\Delta t := \Delta t_1 = \Delta t_2$. We will see later in how far the analysis performed for the nonmultirate case is applicable to the multirate case.

**Algorithm 6** Solver for the Neumann problems in (28).
---
1: **procedure** SDIRK2NEUMANN($\psi_I^{k+1,n}$, $\psi_\Gamma^{k+1,n}$, $\mathbf{F}_j^{k,n}$)
2:     **for** $j = 1, 2$ **do**   *# loop over stages*
3:         $\mathbf{s}_j \leftarrow \psi^{k+1,n} + \Delta t \sum_{l=1}^{j-1}(1-a)\mathbf{k}_l$
4:         $\mathbf{Y}_j \leftarrow$ SOLVELINEARSYSTEM($\mathbf{s}_j, \psi_I^{k+1,n}, \psi_\Gamma^{k+1,n}, \mathbf{F}_j^{k,n}$),   *# solve 1st eq in (28)*
5:         $\mathbf{k}_j \leftarrow \frac{1}{a\Delta t}(\mathbf{Y}_j - \mathbf{s}_j)$
6:     $\psi_I^{k+1,n+1}, \psi_\Gamma^{k+1,n+1} \leftarrow \mathbf{Y}_2$
        **return** $\psi_I^{k+1,n+1}, \psi_\Gamma^{k+1,n+1}$
---

**Algorithm 7** NNWR algorithm using SDIRK2.
---
1: **procedure** NNWR2($\tau_1, \tau_2, \alpha_1, \alpha_2, \lambda_1, \lambda_2, \Theta, TOL$)
2:     $\mathbf{u}_I^{(m),k+1,0}, \mathbf{u}_\Gamma^0(\tau_m), \psi_I^{(m),k+1,n_m+1}, \psi_\Gamma^{(m),k+1,n_m+1} \leftarrow$ INITIALIZATION
3:     **while** $\|\mathbf{u}_\Gamma^{k+1,N_m} - \mathbf{u}_\Gamma^{k,N_m}\| \leq TOL$ **do**
4:         **for** $j = 1, 2$ **do**   *# loop over stages*
5:             **for** $t_{n_1} \in \tau_1$ **do**
6:                 $\mathbf{u}_I^{(1),k+1,n_1+1}, \mathbf{f}_j^{(1),k,n_1} \leftarrow$ SDIRK2DIRICHLET($\mathbf{u}_I^{(1),k+1,n_1}, \mathbf{u}_\Gamma^{k,n_1}, \mathbf{u}_\Gamma^{k,n_1+1}$)
7:             **for** $t_{n_2} \in \tau_2$ **do** (in parallel to 5)
8:                 $\mathbf{u}_I^{(2),k+1,n_2+1}, \mathbf{f}_j^{(2),k,n_2} \leftarrow$ SDIRK2DIRICHLET($\mathbf{u}_I^{(2),k+1,n_2}, \mathbf{u}_\Gamma^{k,n_2}, \mathbf{u}_\Gamma^{k,n_2+1}$)
9:             $\mathbf{F}_j^{(1),k,\tau_1} \leftarrow \mathbf{f}_j^{(1),k,\tau_1} + \mathrm{I}(\tau_1, \tau_2, \mathbf{f}_j^{(2),k,\tau_2})$
10:            $\mathbf{F}_j^{(2),k,\tau_2} \leftarrow \mathbf{f}_j^{(2),k,\tau_2} + \mathrm{I}(\tau_2, \tau_1, \mathbf{f}_j^{(1),k,\tau_1})$  (in parallel to 9)
11:            **for** $t_{n_1} \in \tau_1$ **do**
12:                $\psi_I^{k+1,n_1+1}, \psi_\Gamma^{k+1,n_1+1} \leftarrow$ SDIRK2NEUMANN($\psi_I^{k+1,n_1}, \psi_\Gamma^{k+1,n_1}, \mathbf{F}_j^{(1),k,n_1}$)
13:            **for** $t_{n_2} \in \tau_2$ **do** (in parallel to 11)
14:               $\psi_I^{k+1,n_2+1}, \psi_\Gamma^{k+1,n_2+1} \leftarrow$ SDIRK2NEUMANN($\psi_I^{k+1,n_2}, \psi_\Gamma^{k+1,n_2}, \mathbf{F}_j^{(2),k,n_2}$)
15:         $\mathbf{u}_\Gamma^{k+1,\tau_1} \leftarrow \mathbf{u}_\Gamma^{k,\tau_1} - \Theta\left(\psi_\Gamma^{(1),k+1,\tau_1} + \mathrm{I}(\tau_1, \tau_2, \psi_\Gamma^{(2),k+1,\tau_2})\right)$
16:         $\mathbf{u}_\Gamma^{k+1,\tau_2} \leftarrow \mathbf{u}_\Gamma^{k,\tau_2} - \Theta\left(\psi_\Gamma^{(2),k+1,\tau_2} + \mathrm{I}(\tau_2, \tau_1, \psi_\Gamma^{(1),k+1,\tau_1})\right)$ (p. to 15)
---

The goal now is to find the iteration matrix $\Sigma$ with respect to the final synchronization point $\mathbf{u}_\Gamma^{N_m} \approx \mathbf{u}_\Gamma(T_f)$ because the global error over the time window $[T_0, T_f]$ is assumed to be increasing, having its maximum at the final time $T_f$. Thus, we will find $\Sigma$ such that

$$\mathbf{u}_\Gamma^{k+1,N_m} = \Sigma \mathbf{u}_\Gamma^{k,N_m} + \sum_{i=1}^{2} \left( \varphi^{k+1,\tilde{\tau}_i} + \varphi^{k,\tau_i} \right), \tag{33}$$

where $\varphi^{k,\tau_m}$ are terms dependent on solutions at the previous fixed point iteration $k$ for the time grids $\tau_m = \{t_1, t_2, .., t_{N_m}\}$, $m = 1, 2$ and $\varphi^{k+1,\tilde{\tau}_m}$ are terms dependent on solutions at the current iteration $k+1$ but for the time grids $\tilde{\tau}_m = \{t_1, t_2, .., t_{N_m-1}\} \subset \tau_m$, $m = 1, 2$. To perform the analysis, we neglect all the solutions at previous time steps (indicated by $\varphi^{k+1,\tilde{\tau}_m}$). Thus, we do not to find the exact rate of convergence when having more than one single time step, but instead a good estimate.

We now rewrite (18), (20) and (21)-(22) as an iteration for $\mathbf{u}_\Gamma^{k+1,N_m}$. As we chose above, we omit all the terms in (33) except for the first two. We isolate the term $\mathbf{u}_I^{(m),k+1,N_m}$ from (18) and $\psi_I^{(m),k+1,N_m}$ from the first equation in (20) leading to

$$\mathbf{u}_I^{(m),k+1,N_m} = -\left( \frac{\mathbf{M}_{II}^{(m)}}{\Delta t} + \mathbf{A}_{II}^{(m)} \right)^{-1} \left( \frac{\mathbf{M}_{I\Gamma}^{(m)}}{\Delta t} + \mathbf{A}_{I\Gamma}^{(m)} \right) \mathbf{u}_\Gamma^{k,N_m}, \tag{34}$$

$$\psi_I^{(m),k+1,N_m} = -\left( \frac{\mathbf{M}_{II}^{(m)}}{\Delta t} + \mathbf{A}_{II}^{(m)} \right)^{-1} \left( \frac{\mathbf{M}_{I\Gamma}^{(m)}}{\Delta t} + \mathbf{A}_{I\Gamma}^{(m)} \right) \psi_\Gamma^{(m),k+1,N_m}. \tag{35}$$

Inserting (34) and (35) into the second equation of (20) we get

$$\psi_\Gamma^{(m),k+1,N_m} = \mathbf{S}^{(m)^{-1}} \sum_{i=1}^{2} \mathbf{S}^{(i)} \mathbf{u}_\Gamma^{k,N_m}, \tag{36}$$

with

$$\mathbf{S}^{(m)} := \left( \frac{\mathbf{M}_{\Gamma\Gamma}^{(m)}}{\Delta t} + \mathbf{A}_{\Gamma\Gamma}^{(m)} \right) - \left( \frac{\mathbf{M}_{\Gamma I}^{(m)}}{\Delta t} + \mathbf{A}_{\Gamma I}^{(m)} \right) \left( \frac{\mathbf{M}_{II}^{(m)}}{\Delta t} + \mathbf{A}_{II}^{(m)} \right)^{-1} \left( \frac{\mathbf{M}_{I\Gamma}^{(m)}}{\Delta t} + \mathbf{A}_{I\Gamma}^{(m)} \right). \tag{37}$$

Finally, inserting (36) into (21) or (22) one gets $\mathbf{u}_\Gamma^{k+1,N_m} = \Sigma \mathbf{u}_\Gamma^{k,N_m}$ with

$$\Sigma = \mathbf{I} - \Theta \left( 2\mathbf{I} + \mathbf{S}^{(1)^{-1}} \mathbf{S}^{(2)} + \mathbf{S}^{(2)^{-1}} \mathbf{S}^{(1)} \right). \tag{38}$$

In the one-dimensional case, the iteration matrix $\Sigma$ is just a real number and thus its spectral radius is its modulus. Then, the optimal relaxation parameter $\Theta_{opt}$ in 1D is given by

$$\Theta_{opt} = \frac{1}{2 + \mathbf{S}^{(1)^{-1}} \mathbf{S}^{(2)} + \mathbf{S}^{(2)^{-1}} \mathbf{S}^{(1)}}. \tag{39}$$

## 9 One-dimensional convergence analysis

So far, the derivation was performed for a rather general discretization. In this section, we study the iteration matrix $\Sigma$ for a specific FE discretization in 1D. We will give a formula for the convergence rates. The behaviour of the rates when approaching both the continuous case in time and space is also given.

Specifically, we use $\Omega_1 = [-1, 0]$, $\Omega_2 = [0, 1]$. For the FE discretization, we use the standard piecewise-linear polynomials as test functions. Here we discretize $\Omega_m$ into $N + 1$ equal sized cells of size $\Delta x = 1/(N+1)$ for $m = 1, 2$.

With $\mathbf{e}_j = \begin{pmatrix} 0 & \cdots & 0 & 1 & 0 & \cdots & 0 \end{pmatrix}^T \in \mathbb{R}^N$ where the only nonzero entry is located at the $j$-th position, the discretization matrices are given by

$$\mathbf{A}_{II}^{(m)} = \frac{\lambda_m}{\Delta x^2} \begin{pmatrix} 2 & -1 & & 0 \\ -1 & 2 & \ddots & \\ & \ddots & \ddots & -1 \\ 0 & & -1 & 2 \end{pmatrix}, \quad \mathbf{M}_{II}^{(m)} = \frac{\alpha_m}{6} \begin{pmatrix} 4 & 1 & & 0 \\ 1 & 4 & \ddots & \\ & \ddots & \ddots & 1 \\ 0 & & 1 & 4 \end{pmatrix},$$

$$\mathbf{M}_{\Gamma\Gamma}^{(m)} = \frac{2\alpha_m}{6}, \quad \mathbf{A}_{\Gamma\Gamma}^{(m)} = \frac{\lambda_m}{\Delta x^2},$$

$$\mathbf{A}_{I\Gamma}^{(1)} = -\frac{\lambda_1}{\Delta x^2} \mathbf{e}_N, \quad \mathbf{A}_{I\Gamma}^{(2)} = -\frac{\lambda_2}{\Delta x^2} \mathbf{e}_1, \quad \mathbf{M}_{I\Gamma}^{(1)} = \frac{\alpha_1}{6} \mathbf{e}_N, \quad \mathbf{M}_{I\Gamma}^{(2)} = \frac{\alpha_2}{6} \mathbf{e}_1,$$

$$\mathbf{M}_{\Gamma I}^{(1)} = \frac{\alpha_1}{6} \mathbf{e}_N^T, \quad \mathbf{M}_{\Gamma I}^{(2)} = \frac{\alpha_2}{6} \mathbf{e}_1^T, \quad \mathbf{A}_{\Gamma I}^{(1)} = -\frac{\lambda_1}{\Delta x^2} \mathbf{e}_N^T, \quad \mathbf{A}_{\Gamma I}^{(2)} = -\frac{\lambda_2}{\Delta x^2} \mathbf{e}_1^T.$$

where $\mathbf{A}_{II}^{(m)}, \mathbf{M}_{II}^{(m)} \in \mathbb{R}^{N \times N}$, $\mathbf{A}_{I\Gamma}^{(m)}, \mathbf{M}_{I\Gamma}^{(m)} \in \mathbb{R}^{N \times 1}$ and $\mathbf{A}_{\Gamma I}^{(m)}, \mathbf{M}_{\Gamma I}^{(m)} \in \mathbb{R}^{1 \times N}$ for $m = 1, 2$.

One computes $\mathbf{S}^{(m)}$ for $m = 1, 2$, by inserting the corresponding matrices specified above in (37) obtaining

$$\mathbf{S}^{(1)} = \left(\frac{\alpha_1}{3\Delta t} + \frac{\lambda_1}{\Delta x^2}\right) - \left(\frac{\alpha_1}{6\Delta t} - \frac{\lambda_1}{\Delta x^2}\right)^2 \mathbf{e}_N^T \left(\frac{\mathbf{M}_{II}^{(1)}}{\Delta t} + \mathbf{A}_{II}^{(1)}\right)^{-1} \mathbf{e}_N$$
$$= \left(\frac{\alpha_1}{3\Delta t} + \frac{\lambda_1}{\Delta x^2}\right) - \left(\frac{\alpha_1}{6\Delta t} - \frac{\lambda_1}{\Delta x^2}\right)^2 \alpha_{NN}^1, \quad (40)$$

$$\mathbf{S}^{(2)} = \left(\frac{\alpha_2}{3\Delta t} + \frac{\lambda_2}{\Delta x^2}\right) - \left(\frac{\alpha_2}{6\Delta t} - \frac{\lambda_2}{\Delta x^2}\right)^2 \mathbf{e}_1^T \left(\frac{\mathbf{M}_{II}^{(2)}}{\Delta t} + \mathbf{A}_{II}^{(2)}\right)^{-1} \mathbf{e}_1$$
$$= \left(\frac{\alpha_2}{3\Delta t} + \frac{\lambda_2}{\Delta x^2}\right) - \left(\frac{\alpha_2}{6\Delta t} - \frac{\lambda_2}{\Delta x^2}\right)^2 \alpha_{11}^2, \quad (41)$$

where $\alpha_{ij}^m$ represent the entries of the matrices $\left(\frac{\mathbf{M}_{II}^{(m)}}{\Delta t} + \mathbf{A}_{II}^{(m)}\right)^{-1}$ for $i, j = 1, ..., N$, $m = 1, 2$. Observe that the matrices $\left(\frac{\mathbf{M}_{II}^{(1)}}{\Delta t} + \mathbf{A}_{II}^{(1)}\right)$ and $\left(\frac{\mathbf{M}_{II}^{(2)}}{\Delta t} + \mathbf{A}_{II}^{(2)}\right)$ are tridiagonal Toeplitz matrices but their inverses are full matrices. The computation of the exact inverses could be performed based on the recursive formula presented in [9] which runs over the entries of the matrix and consequently, it is non trivial to compute $\alpha_{NN}^1$ and $\alpha_{11}^2$ this way.

Due to these difficulties, we rewrite them in terms of their eigendecomposition:

$$\left(\frac{\mathbf{M}_{II}^{(m)}}{\Delta t} + \mathbf{A}_{II}^{(m)}\right)^{-1} = \left[\text{tridiag}\left(\frac{\alpha_m \Delta x^2 - 6\lambda_m \Delta t}{6\Delta t \Delta x^2}, \frac{2\alpha_m \Delta x^2 + 6\lambda_m \Delta t}{3\Delta t \Delta x^2}, \frac{\alpha_m \Delta x^2 - 6\lambda_m \Delta t}{6\Delta t \Delta x^2}\right)\right]^{-1}$$
$$= \mathbf{V} \Lambda_m^{-1} \mathbf{V}, \quad \text{for} \quad m = 1, 2,$$

where the matrix $\mathbf{V}$ has the eigenvectors of any symmetric tridiagonal Toeplitz matrix of dimension $N$ as columns. The entries of $\mathbf{V}$ are not dependent on the entries of $\frac{\mathbf{M}_{II}^{(m)}}{\Delta t} + \mathbf{A}_{II}^{(m)}$ due to their symmetry. Moreover, the matrices $\Lambda_m$ are diagonal matrices having the eigenvalues of $\frac{\mathbf{M}_{II}^{(m)}}{\Delta t} + \mathbf{A}_{II}^{(m)}$ as entries. These are known and given e.g. in [25, pp. 514-516]:

$$v_{ij} = \frac{1}{\sqrt{\sum_{k=1}^N \sin^2\left(\frac{k\pi}{N+1}\right)}} \sin\left(\frac{ij\pi}{N+1}\right),$$

$$\lambda_j^m = \frac{1}{3\Delta t \Delta x^2}\left(2\alpha_m \Delta x^2 + 6\lambda_m \Delta t + (\alpha_m \Delta x^2 - 6\lambda_m \Delta t) \cos\left(\frac{j\pi}{N+1}\right)\right),$$

for $i, j = 1, ..., N$ and $m = 1, 2$.

The entries $\alpha^1_{NN}$ and $\alpha^2_{11}$ of the matrices $\left(\frac{\mathbf{M}^{(1)}_{II}}{\Delta t} + \mathbf{A}^{(1)}_{II}\right)^{-1}$ and $\left(\frac{\mathbf{M}^{(2)}_{II}}{\Delta t} + \mathbf{A}^{(2)}_{II}\right)^{-1}$, respectively, are now computed through their eigendecomposition resulting in

$$\alpha^1_{NN} = \frac{\sum_{i=1}^N \frac{1}{\lambda^1_i}\sin^2\left(\frac{i\pi N}{N+1}\right)}{\sum_{i=1}^N \sin^2\left(\frac{i\pi}{N+1}\right)} = \frac{s_1}{\sum_{i=1}^N \sin^2(i\pi\Delta x)}, \qquad (42)$$

$$\alpha^2_{11} = \frac{\sum_{i=1}^N \frac{1}{\lambda^2_i}\sin^2\left(\frac{i\pi}{N+1}\right)}{\sum_{i=1}^N \sin^2\left(\frac{i\pi}{N+1}\right)} = \frac{s_2}{\sum_{i=1}^N \sin^2(i\pi\Delta x)}, \qquad (43)$$

with

$$s_m = \sum_{i=1}^N \frac{3\Delta t \Delta x^2 \sin^2(i\pi\Delta x)}{2\alpha_m \Delta x^2 + 6\lambda_m \Delta t + (\alpha_m \Delta x^2 - 6\lambda_m \Delta t)\cos(i\pi\Delta x)}, \qquad (44)$$

for $m = 1, 2$.

To simplify this, the finite sum $\sum_{i=1}^N \sin^2(i\pi\Delta x)$ can be computed. We first rewrite the sum of squared sinus terms into a sum of cosinus terms using the identity $\sin^2(x/2) = (1-\cos(x))/2$. Then, the resulting sum can be converted into a geometric sum using Euler's formula:

$$\sum_{i=1}^N \sin^2(i\pi\Delta x) = \frac{1-\Delta x}{2\Delta x} - \frac{1}{2}\sum_{i=1}^N \cos(2i\pi\Delta x) = \frac{1}{2\Delta x}. \qquad (45)$$

Now, inserting (45) into (42) and (43) and these two into (40) and (41) we get for $\mathbf{S}^{(m)}$ for $m = 1, 2$,

$$\mathbf{S}^{(m)} = \frac{6\Delta t \Delta x(\alpha_m \Delta x^2 + 3\lambda_m \Delta t) - (\alpha_m \Delta x^2 - 6\lambda_m \Delta t)^2 s_m}{18\Delta t^2 \Delta x^3}, \qquad (46)$$

With this we obtain an explicit formula for the optimal relaxation parameter $\Theta_{opt}$ in (39):

$$\Theta_{opt} = \left(2 + \frac{6\Delta t \Delta x(\alpha_2 \Delta x^2 + 3\lambda_2 \Delta t) - (\alpha_2 \Delta x^2 - 6\lambda_2 \Delta t)^2 s_2}{6\Delta t \Delta x(\alpha_1 \Delta x^2 + 3\lambda_1 \Delta t) - (\alpha_1 \Delta x^2 - 6\lambda_1 \Delta t)^2 s_1}\right.$$
$$\left.+ \frac{6\Delta t \Delta x(\alpha_1 \Delta x^2 + 3\lambda_1 \Delta t) - (\alpha_1 \Delta x^2 - 6\lambda_1 \Delta t)^2 s_1}{6\Delta t \Delta x(\alpha_2 \Delta x^2 + 3\lambda_2 \Delta t) - (\alpha_2 \Delta x^2 - 6\lambda_2 \Delta t)^2 s_2}\right)^{-1}. \qquad (47)$$

We could not find a way of simplifying the finite sum (44) because $\Delta x$ depends on $N$ (i.e., $\Delta x = 1/(N+1)$). However, (47) is a computable expression that gives the optimal relaxation parameter $\Theta_{opt}$ of the NNWR algorithm using implicit Euler for given $\Delta x$, $\Delta t$, $\alpha_m$ and $\lambda_m$, $m = 1, 2$.

We are now interested in the asymptotics of (47) with respect to both spatial and temporal resolutions. To this end, we rewrite (47) in terms of $c := \Delta t / \Delta x^2$:

$$\Theta_{opt} = \left(2 + \frac{6\Delta t(\alpha_2 + 3\lambda_2 c) - \Delta x(\alpha_2 - 6\lambda_2 c)^2 s_2'}{6\Delta t(\alpha_1 + 3\lambda_1 c) - \Delta x(\alpha_1 - 6\lambda_1 c)^2 s_1'} \right. \\ \left. + \frac{6\Delta t(\alpha_1 + 3\lambda_1 c) - \Delta x(\alpha_1 - 6\lambda_1 c)^2 s_1'}{6\Delta t(\alpha_2 + 3\lambda_2 c) - \Delta x(\alpha_2 - 6\lambda_2 c)^2 s_2'}\right)^{-1}. \tag{48}$$

where

$$s_m' = \sum_{i=1}^{N} \frac{3\Delta t \sin^2(i\pi \Delta x)}{2\alpha_m + 6\lambda_m c + (\alpha_m - 6\lambda_m c)\cos(i\pi \Delta x)}, \tag{49}$$

for $m = 1, 2$.

Finally, we compute the limits $c \to 0$ and $c \to \infty$ of the iteration matrix $\Sigma$:

$$\lim_{c \to 0} \Sigma = \lim_{c \to 0} \left(1 - \Theta \left(2 + \frac{6\alpha_2 \Delta t - \alpha_2 \Delta t \Delta x \sum_{i=1}^{N} \frac{3\sin(i\pi\Delta x)^2}{2+\cos(i\pi\Delta x)}}{6\alpha_1 \Delta t - \alpha_1 \Delta t \Delta x \sum_{i=1}^{N} \frac{3\sin(i\pi\Delta x)^2}{2+\cos(i\pi\Delta x)}} \right.\right. \\ \left.\left. + \frac{6\alpha_1 \Delta t - \alpha_1 \Delta t \Delta x \sum_{i=1}^{N} \frac{3\sin(i\pi\Delta x)^2}{2+\cos(i\pi\Delta x)}}{6\alpha_2 \Delta t - \alpha_2 \Delta t \Delta x \sum_{i=1}^{N} \frac{3\sin(i\pi\Delta x)^2}{2+\cos(i\pi\Delta x)}}\right)\right) \\ = 1 - \Theta\left(2 + \frac{\alpha_2}{\alpha_1} + \frac{\alpha_1}{\alpha_2}\right) = 1 - \Theta\left(\frac{(\alpha_1 + \alpha_2)^2}{\alpha_1 \alpha_2}\right),$$

and consequently

$$\Theta_{\{c_1 \to 0\}} = \frac{\alpha_1 \alpha_2}{(\alpha_1 + \alpha_2)^2}. \tag{50}$$

$$\lim_{c \to \infty} \Sigma = \lim_{c \to \infty} \left(1 - \Theta \left(2 + \frac{18\lambda_2 \Delta t c - 18\lambda_2 \Delta t \Delta x c \sum_{i=1}^{N} \frac{\sin(i\pi\Delta x)^2}{1-\cos(i\pi\Delta x)}}{18\lambda_1 \Delta t c - 18\lambda_1 \Delta t \Delta x c \sum_{i=1}^{N} \frac{\sin(i\pi\Delta x)^2}{1-\cos(i\pi\Delta x)}} \right.\right. \\ \left.\left. + \frac{18\lambda_1 \Delta t c - 18\lambda_1 \Delta t \Delta x c \sum_{i=1}^{N} \frac{\sin(i\pi\Delta x)^2}{1-\cos(i\pi\Delta x)}}{18\lambda_2 \Delta t c - 18\lambda_2 \Delta t \Delta x c \sum_{i=1}^{N} \frac{\sin(i\pi\Delta x)^2}{1-\cos(i\pi\Delta x)}}\right)\right) \\ = 1 - \Theta\left(2 + \frac{\lambda_2}{\lambda_1} + \frac{\lambda_1}{\lambda_2}\right) = 1 - \Theta\left(\frac{(\lambda_1 + \lambda_2)^2}{\lambda_1 \lambda_2}\right),$$

and consequently

$$\Theta_{\{c_1 \to \infty\}} = \frac{\lambda_1 \lambda_2}{(\lambda_1 + \lambda_2)^2}. \tag{51}$$

The result obtained in (51) is consistent with the one-dimensional semidiscrete analysis performed in [21]. There, a convergence analysis for the NNWR method in (6), (7) and (8) with constant coefficients using Laplace transforms shows that $\Theta_{opt} = 1/4$ when the two subdomains $\Omega_1$ and $\Omega_2$ are identical. Their result is recovered by our analysis when one approaches the continuous case in space in (51) for constant coefficients, i.e, $\lambda_1 = \lambda_2$. In that case, one gets

$$\Theta_{opt} = \frac{\lambda_1 \lambda_2}{(\lambda_1 + \lambda_2)^2} = \frac{\lambda_1^2}{4\lambda_1^2} = \frac{1}{4}. \tag{52}$$

## 10 Numerical results

We now present numerical experiments to illustrate the validity of the theoretical results of the previous sections. All the results in this section have been produced by implementing algorithms 4 and 7 in Python using the FE discretization specified in the previous section. Firstly, we will show numerically that the NNWR algorithm using implicit Euler preserves first order and using SDIRK2 second order. Secondly, we will show the validity of (47) as estimator for the optimal relaxation parameter $\Theta_{opt}$ of the NNWR algorithm using implicit Euler. We will also show that (47) is a good estimator for the multirate case both using implicit Euler and SDIRK2 and also for 2D examples. Furthermore, we will also show that the theoretical asymptotics deduced both in (50) and (51) match with the numerical experiments. Finally, a comparison between the Dirichlet-Neumann and the Neuman-Neumann methods is included.

### 10.1 NNWR results

Figure 4 shows the error plots of the NNWR algorithm for the coupling of different materials using both implicit Euler and SDIRK2. Physical properties of the materials are shown in table 1. We have considered two initial time grids (for $\Delta t_c = 1/10$ and $\Delta t_f = 1/100$ given), which we then refine several times by a factor of 2:

- (C-C): Two coarse conforming time grids with $\Delta t_1 = \Delta t_2 = \Delta t_c$.

- (C-F): Nonconforming time grids with $\Delta t_1 = \Delta t_c$ and $\Delta t_2 = \Delta t_f$.

When coupling two different materials in the multirate case, we always assign the finer grid to the material that has higher heat conductivity because it performs the heat changes faster. In space, we fix $\Delta x = 1/200$ and we compute a reference solution by solving problem (1) directly on a very fine time grid, with $\Delta t = 1/1000$. One observes in figure 4 that first and second order convergence is obtained in the nonconforming case for implicit Euler and SDIRK2 respectively. Moreover, the errors obtained in the multirate case (C-F) are nearly the same as in the coarser nonmultirate case (C-C). Thus, the accuracy of the multirate case is determined by its coarser rate. This is consistent with [31, 8] where the convergence of the discrete multirate WR algorithm is independent of the ratio of timesteps.

Table 1: Physical properties of the materials. $\lambda$ is the thermal conductivity, $\rho$ the density, $c_p$ the specific heat capacity and $\alpha = \rho c_p$.

| Material | $\lambda$ (W/mK) | $\rho$ (kg/m$^3$) | $c_p$ (J/kgK) | $\alpha$ (J/K m$^3$) |
|:---:|:---:|:---:|:---:|:---:|
| **Air** | 0.0243 | 1.293 | 1005 | 1299.5 |
| **Water** | 0.58 | 999.7 | 4192.1 | 4.1908e6 |
| **Steel** | 48.9 | 7836 | 443 | 3471348 |

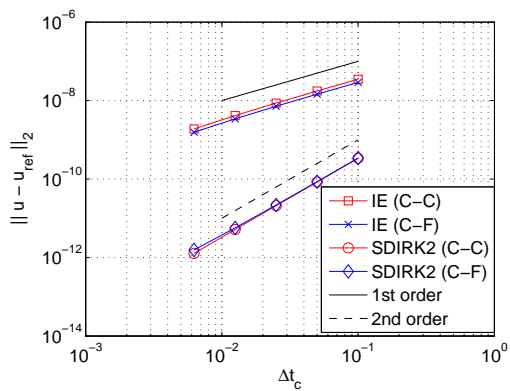

Figure 4: Error plot of the NNWR algorithm for three different time grids. $\Delta x = 1/200$, $[T_0, T_f] = [0, 1]$ and $TOL = 1e - 15$. Air-steel thermal interaction using implicit Euler and air-water thermal interaction using SDIRK2.

Figure 5 compares the behavior of the algorithm described in this paper using implicit Euler (left plot) and SDIRK2 (right plot). It shows the convergence rates in terms of the relaxation parameter $\Theta$ for the one-dimensional thermal coupling between air and water. We have plotted $\Sigma(\Theta)$ in (38) with the 1D space discretization specified in section 9 and the experimental convergence rates for both the multirate and nonmultirate cases. The relevance of the analysis presented above is illustrated in figure 5 because the algorithm is extremely fast at $\Theta_{opt}$ (converging in 2 iterations), but if one deviates slightly from $\Theta_{opt}$, we get a divergent method. As can be seen in the left plot, the experimental convergence rates for the nonmultirate case (C-C) are exactly predicted by the theory. Moreover, our formula also predicts where the convergence rate of the NNWR algorithm in the multirate case (C-F) is minimal. They are not identical because the linear interpolation performed at the space-time interface in the multirate case is neglected in (33) to simplify the theoretical analysis. One can also observe in the right plot that $\Sigma(\Theta)$ using implicit Euler estimates quite well the optimal relaxation parameter of the NNWR algorithm using SDIRK2 for both the multirate and nonmultirate cases.

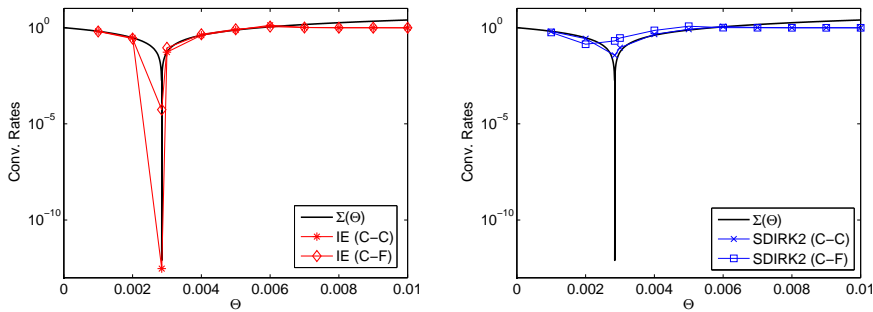

Figure 5: Air-water convergence rates as a function of the relaxation parameter $\Theta$ in 1D. $\Delta x = 1/100$, $\Delta t_c = 100$ and $\Delta t_f = 1$. *Left*: $\Sigma(\Theta)$ in (38) and the experimental convergence rates both for the multirate (C-F) and nonmultirate (C-C) cases using implicit Euler. *Right*: $\Sigma(\Theta)$ in (38) using implicit Euler and the experimental convergence rates both for the multirate (C-F) and nonmultirate (C-C) cases using SDIRK2.

In order to illustrate the behavior of the NNWR method in the multirate case ($\Delta t_1 \neq \Delta t_2$), we have plotted in figure 6 the convergence rates using the relaxation parameters $\Theta_{opt}(\Delta t_1)$ and $\Theta_{opt}(\Delta t_2)$ in (47) with respect to the variation of $\Delta t_1/\Delta t_2$ for the air-water coupling. In 6 we have chosen $\Delta t_1/\Delta t_2 = 1e-3/2e-1, 2e-3/2e-1, 1e-2/2e-1, 2e-2/2e-1, 5e-2/2e-1, 1e-1/2e-1, 2e-1/2e-1, 2e-1/1e-1, 2e-1/5e-2, 2e-1/2e-2, 2e-1/1e-2, 2e-1/2e-3, 2e-1/1e-3$ and $\Delta x = 1/100$. One observes that the multirate method converges fast using any of the two relaxation parameters for both implicit Euler and SDIRK2. Nevertheless, one can also see in figure 6 that even though we have not performed a specific analysis for the optimal relaxation parameter in the multirate case, the $\Theta_{opt}$ in (47) can be used as an estimator. More specifically, we conclude from figure 6 that one can use $\Theta_{opt}(\Delta t_2)$ when $\Delta t_1 < \Delta t_2$ and $\Theta_{opt}(\Delta t_1)$ when $\Delta t_1 > \Delta t_2$.

Figure 7 shows the optimal relaxation parameter $\Theta_{opt}$ with respect to the parameter $c := \Delta t/\Delta x^2$ using both implicit Euler and SDIRK2. We have chosen $\Delta t = 1e-9, 1e-8, ..., 1e8, 1e9$ and $\Delta x = 1/100$. For implicit Euler, we have plotted the function $\Theta_{opt}(c)$ in (48). For SDIRK2, we have plotted a sister function $\Theta_{opt}(c)$ that can be found applying exactly the derivation presented in sections 8 and 9 to the discretized SDIRK2-NNWR method introduced in sec-tion 7.2. One can see that the two time discretization methods have a similar behavior when varying $c$. This illustrates why the optimal relaxation parameter $\Theta_{opt}$ computed in (47) for implicit Euler is also valid for SDIRK2 as observed in figure 5. Furthermore, in 7 we observe that the optimal relaxation parameter for any given $\Delta t$ and $\Delta x$ is always between the bounds of the theoretical asymptotics deduced both in (50) and (51), tending to them in the temporal and spatial limits respectively.

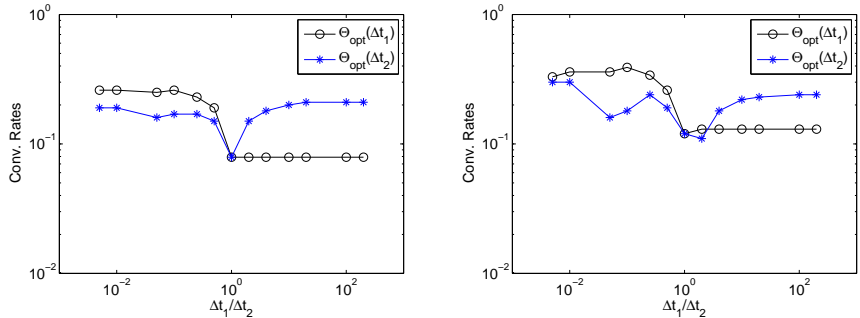

Figure 6: Convergence rates as a function of the temporal ratio $\Delta t_1/\Delta t_2$ for the air-water coupling in 1D. We plot the convergence rates in the multirate case ($\Delta t_1 \neq \Delta t_2$) using the relaxation parameters $\Theta_{opt}(\Delta t_1)$ and $\Theta_{opt}(\Delta t_2)$ in (47). $\Delta t_1/\Delta t_2 = 1e-3/2e-1, 2e-3/2e-1, 1e-2/2e-1, 2e-2/2e-1, 5e-2/2e-1, 1e-1/2e-1, 2e-1/2e-1, 2e-1/1e-1, 2e-1/5e-2, 2e-1/2e-2, 2e-1/1e-2, 2e-1/2e-3, 2e-1/1e-3$ and $\Delta x = 1/100$. *Left*: Implicit Euler. *Right*: SDIRK2.

We now want to demonstrate that the 1D formula (47) is a good estimator for the optimal relaxation parameter $\Theta_{opt}$ in 2D. Thus, we now consider a 2D version of the model problem consisting of two coupled linear heat equations on two identical unit squares, e.g, $\Omega_1 = [-1, 0] \times [0, 1]$ and $\Omega_2 = [0, 1] \times [0, 1]$.

Figure 8 shows the convergence rates in terms of the relaxation parameter $\Theta$ for 2D examples. On the left we have the thermal coupling between air and steel and on the right between air and water. One can observe that the convergence rates of the NNWR method using $\Theta_{opt}$ from (47) in the 2D examples are worse than in 1D, but still optimal. Hence, we suggest to use $\Theta_{opt}$ in 2D as well, otherwise the method is divergent.

## 10.2  NNWR - DNWR Comparison

Finally, we compare the Dirichlet-Neumann and the Neumann-Neumann couplings. We consider the FE discretization specified in section 9 and the implicit Euler as a time integration method for both DNWR and NNWR with $\Delta x = 1/500$. In addition, we will use $\Theta = 1/2$ as the optimal relaxation parameter for the DNWR algorithm as suggested in [13] for constant material coefficients. Note that the optimal relaxation parameter for the NNWR method is $\Theta_{opt} = 1/4$ (see (47)) in the case of constant coefficients because $\lambda_1 = \lambda_2$ and $\alpha_1 = \alpha_2$.

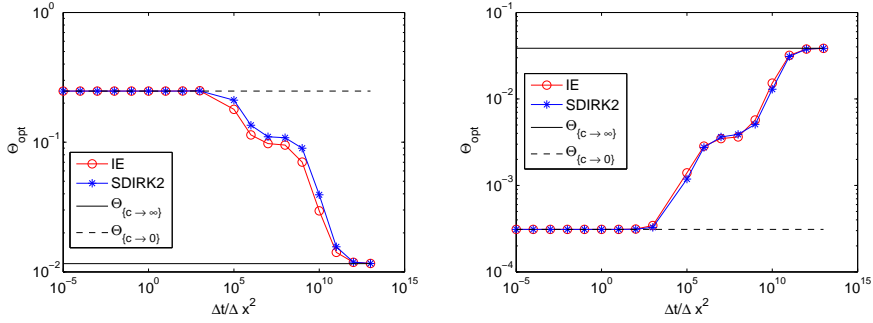

Figure 7: Optimal relaxation parameter $\Theta_{opt}$ as a function of the parameter $c := \Delta t / \Delta x^2$ for both implicit Euler and SDIRK2 in 1D. The constant lines $\Theta_{\{c \to \infty\}}$ and $\Theta_{\{c \to 0\}}$ represent the spatial and temporal asymptotics of $\Theta_{opt}$ in (48). $\Delta t = 1e-9, 1e-8, ..., 1e8, 1e9$ and $\Delta x = 1/100$. *Left*: Water-steel coupling. *Right*: Air-Water coupling.

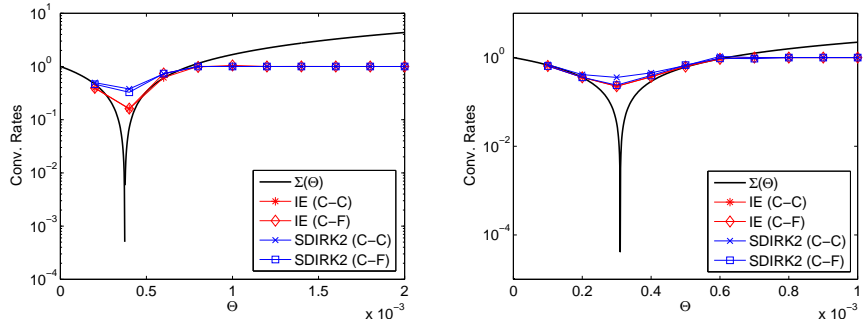

Figure 8: Convergence rates as a function of the relaxation parameter $\Theta$ in 2D. $\Sigma(\Theta)$ in (47) using implicit Euler and the experimental convergence rates both for the multirate (C-F) and the nonmultirate (C-C) cases using implicit Euler and SDIRK2. $\Delta x = 1/10$, $\Delta t_c = 1/10$ and $\Delta t_f = 1/50$. *Left*: Air-steel coupling. *Right*: Air-water coupling.

Table 2 shows the time needed to solve the 1D steel-steel coupling in the nonmultirate case. The number of fixed point iterations needed to achieve a chosen tolerance of $1e-8$ is also given. One can see that the DNWR method is sightly more efficient than the NNWR method. Moreover, the NNWR algorithm runs in parallel on two different processors using double the amount of computational power than the DNWR. Thus, the DNWR method is a better option for this case because it is cheaper and faster.

Table 2: Computational effort of DNWR and NNWR for the 1D steel-steel coupling in the nonmultirate case. $\Delta x = 1/500$ and $TOL = 1e-8$. Number of fixed point iterations in brackets.

| $\Delta t$ | Comp. Time - NNWR (s) | Comp. Time - DNWR (s) |
|---|---|---|
| 1 | 0.07 (2 iterations) | 0.05 (2 iterations) |
| 1/10 | 0.42 (2 iterations) | 0.3 (2 iterations) |
| 1/50 | 1.86 (2 iterations) | 1.25 (2 iterations) |
| 1/100 | 3.74 (2 iterations) | 2.47 (2 iterations) |

However, the NNWR algorithm beats the DNWR algorithm by far when we move to the multirate environment. This is illustrated in table 3 where the computational effort used to solve the 1D steel-steel coupling in the multirate case is shown. There, one can see how the number of fixed point iterations needed to achieve a tolerance of $1e-8$ using DNWR grows exponentially when the difference between $\Delta t_1$ and $\Delta t_2$ increases. On the contrary, the NNWR method is very efficient even when there is a huge difference between $\Delta t_1$ and $\Delta t_2$. Thus, we recommend the NNWR algorithm when coupling two fields with nonconforming time grids.

Table 3: Computational effort comparison of DNWR and NNWR for the 1D steel-steel coupling in the multirate case. $\Delta x = 1/500$ and $TOL = 1e-8$. Number of fixed point iterations in brackets.

| $\Delta t_1 - \Delta t_2$ | Comp. Time - NNWR (s) | Comp. Time - DNWR (s) |
|---|---|---|
| 1/5 - 1/10 | 0.48 (3 iterations) | 0.70 (6 iterations) |
| 1/5 - 1/50 | 1.5 (3 iterations) | 26.98 (71 iterations) |
| 1/5 - 1/100 | 2.74 (3 iterations) | Not convergent |

Finally, we have included a comparison for the 1D air-steel coupling in the multirate case. This interaction between air and steel has the particularity of strong jumps in the material coefficients across the space interface. In this case, we have chosen $\Theta = 1/2$ for DNWR because even though in [13] it is only proved optimal for the case of constant coefficients, they show in the numerical results section that also applies to an example where the diffusion coefficient varies spatially. Moreover, $\Theta_{opt}$ in (47) is chosen for the NNWR method. Table 4 shows a comparison of the computational time employed to solve the 1D air-steel coupling in the multirate case. One can see that the NNWR method is more efficient than the DNWR method because it needs way less fixed point iterations to achieve the same tolerance. Note that the number of iterations increases when the time resolution decreases for the NNWR method. This happens because the analysis performed in section 8 to find the optimal relaxation parameter takes into account only one single time step (see (33)). Therefore, in the case of multiple time steps, $\Theta_{opt}$ in (47) is a very good choice, but it is not optimal. Besides that, the large amount of iterations performed by the DNWR algorithm hints that $\Theta = 1/2$ might not be the optimal relaxation parameter when having strong jumps in the material coefficients for the fully discrete problem. Thus, performing an specific analysis to find the optimal relaxation parameter of the DNWR algorithm is left for future research.

Table 4: Computational effort comparison of DNWR and NNWR for the 1D air-steel coupling in the multirate case. $\Delta x = 1/500$ and $TOL = 1e-8$. Number of fixed point iterations in brackets.

| $\mathbf{\Delta t_1 - \Delta t_2}$ | Comp. Time - NNWR (s) | Comp. Time - DNWR (s) |
|---|---|---|
| 1/5 - 1/10 | 0.47 (3 iterations) | 1.24 (12 iterations) |
| 1/5 - 1/50 | 2.20 (4 iterations) | 4.65 (12 iterations) |
| 1/5 - 1/100 | 3.95 (4 iterations) | 8.77 (12 iterations) |

## 11 Summary and conclusions

We suggested a new high order parallel NNWR method with nonconforming time grids for two heterogeneous coupled heat equations and studied the optimal relaxation parameter in terms of the material coefficients and the temporal and spatial resolutions $\Delta t$ and $\Delta x$. To this end, we considered the coupling of two heat equations on two identical domains. We assumed structured spatial grids and conforming time grids on both subdomains to derive a formula for the optimal relaxation parameter $\Theta_{opt}$ in 1D using implicit Euler. Furthermore, we determined the limits of the optimal relaxation parameter when approaching the continuous case either in space ($\lambda_1\lambda_2/(\lambda_1+\lambda_2)^2$) or time ($\alpha_1\alpha_2/(\alpha_1+\alpha_2)^2$). The method using $\Theta_{opt}$ converges extremely fast, typically within two iterations. This was confirmed by numerical results, where we also demonstrated that the nonmultirate 1D case gives excellent estimates for the multirate 1D case and even for multirate and nonmultirate 2D examples using both implicit Euler and SDIRK2. In addition, we have shown that the NNWR method is a more efficient choice than the classical DNWR in the multirate case.

## References

[1] A. BANKA, *Practical Applications of CFD in heat processing*, Heat Treating Progress., (2005).

[2] P. Birken, T. Gleim, D. Kuhl, and A. Meister, *Fast Solvers for Unsteady Thermal Fluid Structure Interaction*, Int. J. Numer. Meth. Fluids, 79(1) (2015), pp. 16–29.

[3] P. Birken, K. Quint, S. Hartmann, and A. Meister, *A time-adaptive fluid-structure interaction method for thermal coupling*, Comp. Vis. in Science, 13(7) (2011), pp. 331–340.

[4] S. Bremicker-Trübelhorn and S. Ortleb, *On Multirate GARK Schemes with Adaptive Micro Step Sizes for Fluid-Structure Interaction: Order Conditions and Preservation of the Geometric Conservation Law*, Aerospace, 4(8) (2017).

[5] J. Buchlin, *Convective Heat Transfer and Infrared Thermography*, J. Appl. Fluid Mech., 3 (2010), pp. 55–62.

[6] P. Causin, J. Gerbeau, and F. Nobile, *Added-mass effect in the design of partitioned algorithms for fluid-structure problems*, Comp. Meth. Appl. Mech. Eng., 194 (2005), pp. 4506–4527.

[7] A. Cristiano, I. Malossi, P. Blanco, S. Deparis, and A. Quarteroni, *Algorithms for the partitioned solution of weakly coupled fluid models for cardiovascular flows*, Numer. Meth. Biomed. Eng., 27(12) (2011), pp. 2035–2057.

[8] M. Crow and M. Ilic, *The waveform relaxation method for systems of differential/algebraic equations*, in 29th IEEE Conference on Decision and Control, vol. 2, 1990, pp. 453–458.

[9] C. Fonseca and J. Petronilho, *Explicit inverses of some tridiagonal matrices*, Linear Algebra Appl., 325(1-3) (2001), pp. 7–21.

[10] M. Gander, *50 years of time parallel time integration*, in Multiple Shooting and Time Domain Decomposition Methods, T. Carraro, M. Geiger, S. Körkel, and R. Rannacher, eds., Springer, Heidelberg., 2015.

[11] M. Gander, L. Halpern, C. Japhet, and V. Martin, *Advection diffusion problems with pure advection approximation in subregions*, in Domain Decomposition Methods in Science and Engineering XVI. Lecture Notes in Computer Science and Engineering, vol 50, pp. 239-246. Springer, Berlin., 2007.

[12] M. Gander, L. Halpern, and F. Nataf, *Optimized Schwarz waveform relaxation for the one dimensional wave equation*, SIAM J. Numer. Anal., 41(5) (2003), pp. 1643–1681.

[13] M. Gander, F. Kwok, and B. Mandal, *Dirichlet-Neumann and Neumann-Neumann waveform relaxation algorithms for parabolic problems*, ETNA, 45 (2016), pp. 424–456.

[14] M. Gander and A. Stuart, *Space-time continuous analysis of waveform relaxation for the heat equation*, SIAM J. Sci. Comput., 19(6) (1998), pp. 2014–2031.

[15] E. Giladi and H. Keller, *Space-time domain decomposition for parabolic problems*, Numer. Math., 93 (2002), pp. 279–313.

[16] U. Heck, U. Fritsching, and B. K., *Fluid flow and heat transfer in gas jet quenching of a cylinder*, Int. J. Numer. Methods Heat Fluid Flow, 11 (2001), pp. 36–49.

[17] M. Hinderks and R. Radespiel, *Investigation of Hypersonic Gap Flow of a Reentry Nosecap with Consideration of Fluid Structure Interaction*, AIAA Paper, 6 (2006), pp. 2708–3741.

[18] T. Hoang, *Space-time domain decomposition methods for mixed formulations of flow and transport problems in porous media*, PhD thesis, Universiteé Pierre et Marie Curie, Paris 6, France, 2013.

[19] D. Kowollik, P. Horst, and M. Haupt, *Fluid-structure interaction analysis applied to thermal barrier coated cooled rocket thrust chambers with subsequent local investigation of delamination phenomena*, Progress in Propulsion Physics, 4 (2013), pp. 617–636.

[20] D. Kowollik, V. Tini, S. Reese, and M. Haupt, *3D fluid-structure interaction analysis of a typical liquid rocket engine cycle based on a novel viscoplastic damage model*, Int. J. Numer. Methods Engrg., 94 (2013), pp. 1165–1190.

[21] F. Kwok, *Neumann-Neumann waveform relaxation for the time-dependent heat equation*, vol. 98, in Domain Decomposition in Science and Engineering XXI, J. Erhel, M.J. Gander, L. Halpern, G. Pichot, T. Sassi and O.B. Widlund, eds. Lect. Notes Comput. Sci. Eng., 2014.

[22] E. Lelarasmee, A. Ruehli, and A. Sangiovanni-Vincentelli, *The waveform relaxation method for time-domain analysis of large scale integrated circuits*, IEEE Trans. Comput. Aided Des. Integr. Circuits Syst., 1(3) (1982), pp. 131–145.

[23] M. Mehl, B. Uekermann, H. Bijl, D. Blom, B. Gatzhammer, and A. van Zuijlen, *Parallel coupling numerics for partitioned fluid-structure interaction simulations*, Comput. Math. Appl., 71(4) (2016), pp. 869–891.

[24] R. Mehta, *Numerical Computation of Heat Transfer on Reentry Capsules at Mach 5*, AIAA-Paper, 178 (2005).

[25] C. Meyer, *Matrix Analysis and Applied Linear Algebra*, 2000.

[26] A. Monge and P. Birken, *On the convergence rate of the Dirichlet-Neumann iteration for unsteady thermal fluid-structure interaction*, Computational Mechanics, (2017).

[27] A. Quarteroni and A. Valli, *Domain Decomposition Methods for Partial Differential Equations*, Oxford Science Publications, 1999.

[28] P. Stratton, I. Shedletsky, and M. Lee, *Gas Quenching with Helium*, Solid State Phenomena, 118 (2006), pp. 221–226.

[29] A. Toselli and O. Widlund, *Domain Decomposition Methods - Algorithms and Theory*, Springer, 2004.

[30] E. van Brummelen, *Partitioned iterative solution methods for fluid-structure interaction*, International Journal for Numerical Methods in Fluids, 65(1-3) (2011), pp. 3–27.

[31] J. K. White, *Multirate Integration Properties of Waveform Relaxation with Application to Circuit Simulation and Parallel Computation*, PhD thesis, EECS Department, University of California, Berkeley, 1985.